\documentclass[onefignum,onetabnum]{siamltex1213}

\usepackage[utf8]{inputenc}
\usepackage{amsmath,amssymb,amsfonts}

\usepackage{graphicx}
\DeclareGraphicsExtensions{.png,.jpg,.pdf}

\usepackage{xcolor} 
\usepackage{pgfplots}
\pgfplotsset{compat=1.13}

\usepackage{fp}
\def\resrow#1#2#3#4{$#1$&\FPeval{\result}{round(1+(#3-(#4))/#2:5)}\result&$#2$&$#3$&$#4$\\\hline}

\def\R{\mathbb{R}}
\def\d{\partial}
\def\n{\nabla}

\begin{document}

\pagestyle{myheadings}
\slugger{sisc}{xxxx}{xx}{x}{x--x}
\markboth{N.~Rahimi, P.~Kerfriden, F.\,C.~Langbein, R.\,R.~Martin}{CAD Model Simplification Error Estimation for Electrostatics Problems}

\renewcommand{\thefootnote}{\fnsymbol{footnote}}
\title{CAD Model simplification error estimation for electrostatics problems%
  \thanks{This work was supported by the EU FP7 Marie Curie ITN INSIST grant agreement n$^{\small\rm o}$ 289361.}}
\author{N.~Rahimi\footnotemark[2] \and%
        P.~Kerfriden\footnotemark[3] \and%
        F.\,C.~Langbein\footnotemark[2] \and%
        R.\,R.~Martin\footnotemark[2]}%
\footnotetext[2]{School of Computer Science and Informatics, Cardiff University, UK.}
\footnotetext[3]{School of Engineering, Cardiff University, UK.}
\maketitle
\renewcommand{\thefootnote}{\arabic{footnote}}

\begin{abstract}
  Simplifying the geometry of a CAD model using defeaturing techniques enables more efficient discretisation and subsequent simulation for engineering analysis problems. Understanding the effect this simplification has on the solution helps to decide whether the simplification is suitable for a specific simulation problem. It can also help to understand the functional effect of a geometry feature. The effect of the simplification is quantified by a user-defined quantity of interest which is assumed to be (approximately) linear in the solution. A bound on the difference between the quantity of interest of the original and simplified solutions based on the energy norm is derived. The approach is presented in the context of electrostatics problems, but can be applied in general to a range of elliptic partial differential equations. Numerical results on the efficiency of the bound are provided for electrostatics problems with simplifications involving changes inside the problem domain as well as changes to the boundaries.
\end{abstract}


\section{Introduction}

\begin{figure}[t]
 \centering
 \includegraphics[width=\textwidth]{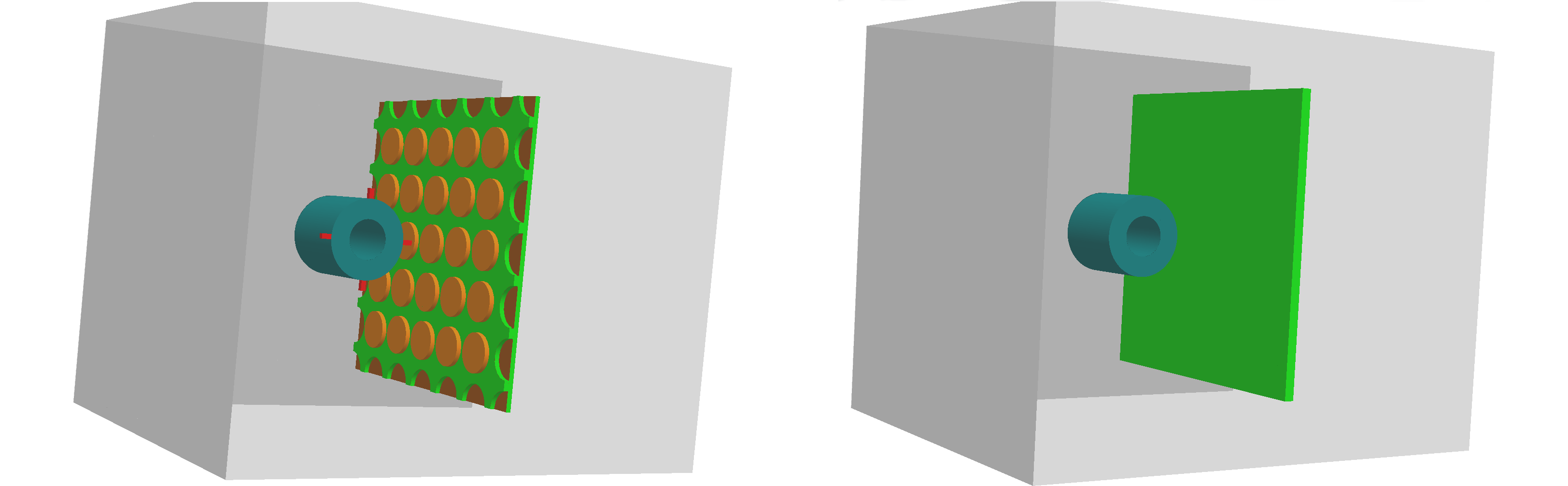}
 \caption{Geometric simplification before magnetostatic analysis of a shielded coil. Left: original fully-featured model.  Right: defeatured model after suppression of unimportant features.}
 \label{fig:3D-Model}
\end{figure}

Computational engineering analysis requires discretization of a continuous boundary value problem. The discretization quality strongly influences the solution accuracy, which depends mainly on (i) how well the properties of the continuous solution space are preserved in the discretized functional solution space and (ii) how well the discrete geometry (typically a 2D or 3D mesh) represents the continuous geometry. It is well known that generating mesh models from CAD models for engineering analysis is time-consuming and expensive, taking $60\%$ to $70\%$ of the total analysis time, because of the algorithms failing to produce a suitable mesh and so manual intervention is required~\cite{IGA-Hughes}. A common approach therefore, is to simplify or idealize the CAD model geometry, removing small or insignificant features which have little effect on the analysis results. This has two advantages: firstly the simpler geometry means that it can be represented by a simpler mesh with fewer, larger elements, making meshing both quicker, and more robust. Secondly, as the resulting mesh is simpler, analysis is also quicker. Fig.~\ref{fig:3D-Model} illustrates an example of simplifying a geometric model of a shielded coil prior to magnetostatic analysis. The model is an example for shielding the magnetic field of a coil. The grey box is the outer boundary which is usually made of metal. The internal components are the coil and shield. The red and orange parts are features on the coil and shield that may be simplified. Much of the geometry has little effect on the solution, and can be removed before meshing.

The key to performing simplification is to know what effect the simplification will have on the solution. Simplifications which make large changes to results are likely to be unacceptable. We wish to predict the simplification error arising from replacing complex geometry by simpler geometry, to decide whether the simplified model is sufficient or a more complex and expensive simulation needs to be run instead. Note that this is problem dependent: features which are important for a strength analysis may be unimportant for electrostatic analysis and vice versa.

Various approaches for simplifying or defeaturing CAD models have been studied~\cite{Danglade,William,Atul}. Their focus is often on the defeaturing techniques themselves, aimed at removing small features, rather than on the effects of defeaturing on subsequent analysis. \cite{Atul}~lists the most common defeaturing techniques applied to industrial designs. There are, however, fewer results on the effect of simplification on the accuracy of the solution. To study this effect, an error measure indicating the effect of the simplification must be defined. This is not simple, as the difference between the solutions for the original and simplified model is not clearly defined, due to differences in geometry and possibly also topology. However, most engineering analysis tasks are aimed at calculating some specific \emph{quantities of interest} (QoIs). Often a QoI can be expressed as an integral of a local quantity, determined by the solution, over a subset of the model, referred to as the domain of interest~\cite{Oden-Prud,Oden}. We call the difference between the QoI of the original and simplified model the \emph{simplification error}. This is the quantity we wish to estimate in order to decide whether the simplified model is suitable for the analysis task. Various examples of industrial applications of defeaturing are collected in~\cite{Shimada}, showing how simplification can significantly improve computational costs in CAE, which makes a strong case for needing good simplification error estimates. Recent results in defining features declaratively with efficient detection algorithms also enable simple access to features specific to the analysis problem~\cite{Niu}.

Estimating the simplification error can be done similarly to determining other error estimates for finite element and boundary element methods. Approximating a real object with a discretized CAD model causes analysis errors; one  aim of error analysis is thus to show that the computational results are close to the exact solution~\cite{Bathe}. To deal with discretization errors, various strategies are proposed for analytical solutions of finite element methods~\cite{Kelly}. \emph{A priori} error estimation studies convergence and stability of numerical solvers, and provides an understanding of the asymptotic behavior of errors for varying mesh parameters~\cite{Oden}. \emph{A posteriori} error estimates in computational mechanics control the error of the solution of ordinary differential equations. Those are typified by predictor-corrector algorithms, where we wish to estimate the differences in errors in solutions obtained by schemes using different orders of error truncation~\cite{Babuska}. The latter approach also provides the basis for our simplification error estimation.

While our approach can be applied to a wide range of physical problems governed by a second order divergence free partial differential equation (PDE), we study it specifically for electrostatics. This provides a simple setting in which we can demonstrate our approach. It is also relevant to the design of capacitors and other problems involving effects of stationary charges, e.g.\ damage to electronic components during manufacturing, and build-up of static electricity, especially on plastics and high-voltage components. Electrostatic analysis is important in the design of switchgear, used to control, protect and isolate electrical equipment. Insulation materials such as gas, oil and air are used; flaws in manufacturing processes may significantly impact its usability. Estimating the manufacturing design error can play a major role in understanding whether the switchgear will operate sufficiently well. Hence, we apply our approach to capacitor analysis problems which can be exemplified in three different schemes:
\begin{enumerate}
\item[(i)] The boundary conditions associated with the problem are changed by simplifying the boundary geometry for practical defeaturing reasons. The change in the boundary conditions results in a different solution.
\item[(ii)] Due to uncertainties in manufacturing processes, the boundaries of a real capacitor are not exactly as intended in the design. Here simplification error estimation predicts the performance of the designed model under realistic, small manufacturing changes in the model geometry.
\item[(iii)] The amount of energy stored in a capacitor depends on the dielectric properties of the material between the capacitor plates and their geometry. Suppressing or modifying the dielectric properties of the material inside the capacitor is sometimes required to correct the dielectric material because of manufacturing uncertainties, and studying the effects on capacitor quality factors or for defeaturing purposes to simplify computation.
\end{enumerate}
We provide guaranteed upper and lower bounds for the simplification error in the energy norm due to such changes for a QoI, which is linear in the solution of the problem. In case the QoI must be linearized for this to hold, there are of course additional approximation errors, for which we assume they are considerably smaller than the simplification error such that they can be ignored. Specific simplifications we consider later relate to changing dielectric material properties in a domain and removing positive or negative features on various boundaries.

Our approach constructs bounds for the simplification error by utilising the concept of the \emph{constitutive relation error} (CRE)~\cite{CRE}. CRE manipulates the admissibility conditions and computes the error measure in the energy norm which makes the computation less expensive. It assumes a QoI which is linear in the solution and provides a link to the residual error in goal-oriented error estimation.

Our simplification error estimation technique and its accuracy are evaluated for capacitor models for schemes (i)-(iii) above with different simplifications. We test it by removing a feature inside the domain, so that the dielectric material constant for an internal region is replaced by the dielectric material constant for the domain surrounding it. This analysis is repeated for several dielectric material setups, feature sizes and locations. Other tests involve removing positive and negative features on the outer boundary or the electronic components (on Neumann or Dirichlet boundaries).

In the rest of the paper, we first review related work and electrostatics boundary value problems. In Section~\ref{General}, our approach toward simplification error estimation is discussed for Laplacian operators, based on goal-oriented error estimation. We show how to derive upper and lower bounds that will consequently be used for various simplification cases. Section~\ref{Inner} employs these bounds for internal features and presents results for various test cases. Section~\ref{Boundary} then shows how to derive simplification error bounds for boundary features with Neumann and Dirichlet boundary conditions and again presents results for test cases. We conclude with a general discussion in Section~\ref{conclusion}.

\section{Related Work}

A posteriori error estimation can be applied to most elliptical equations. It forms an important basis for simplification, modelling and numerical error estimation. For example,~\cite{Ferrandes} analyzes the finite element model for a linear stationary thermal conductivity problem and computes an error estimate for a global quantity of heat conduction related to the solution of the linear heat equation after geometric simplification. In a linear elasticity problem,~\cite{Inna,Inna1} demonstrate how a posteriori error estimation helps in \emph{shape and topological sensitivity analysis} for defeaturing error analysis. Shape sensitivity analysis computes the change in QoIs when the shape of the model is perturbed infinitesimally, while topology sensitivity analysis computes the change when an infinitesimal internal feature is added to or removed from a model. Topological sensitivity provides a powerful technique for shape optimization of arbitrary-shaped features. The final error indicator, however, is estimated only roughly and there is no localization of the error in their technique. Li~\cite{Ming2} extends the work of~\cite{Inna1} to determine the defeaturing error, again in an a posteriori error estimation framework, when defeaturing CAD models for different feature types. They use \emph{shape idealization} for dimension reduction of a thin model to a 2D plate, and calculate the energy norm of the induced error for specific QoIs arising from the solution of a linear Poisson equation. Their margin of error in estimating the simplification error is large and there is no guarantee that the technique provides a strict bound for the actual error.

\cite{Suresh}~describes removing a negative feature for a thermal conductivity equation where the feature is subject to a Neumann boundary condition. This work was extended by~\cite{Ming3} for negative features with Neumann boundary condition. It uses an approach based on \emph{dual weighted residuals} (DWR) which is based on reformulating the modification sensitivity, originally caused by a geometric difference, as a modeling error, caused by mathematical modeling of PDEs over the same geometric model. The DWR is calculated via sensitivity analysis by integrating over the feature's boundary which can be evaluated using engineering analysis results from the defeatured model. This technique is closely related to our method, except for only being able to deal with features on Neumann boundaries.

The above discussed results can only handle the error caused by simplifying a single feature. \cite{Ming6}~shows how to construct a defeaturing error estimator for second-order shape sensitivity that considers the interaction between different internal or boundary features. However, the estimation results are quite inaccurate, and the assumptions made concerning boundary features result in rather wide error estimate ranges. While~\cite{Ming1} uses goal-oriented a posteriori error estimation to improve defeaturing error estimation for Poisson and linear elasticity equations, the approach does not include a strict bound for the error that would prevent over- or underestimation. Thus, in~\cite{Ming4} adjoint theory is employed to estimate the simplification error, which improves the error estimation results for internal and negative boundary features with Neumann and Dirichlet boundary conditions. Green's theorem and linearization are used to derive a simplification error formula for a non-linear Poisson equation. All techniques presented have some drawbacks in that only certain types of feature or boundary conditions can be dealt with, there is often no guarantee that the estimate will bound or be sufficiently close to the actual error, nor is there a consistent index to indicate the performance of the estimate. Also, none of the previous work takes into consideration the solution and physics governed by the electrostatics equation.

Our approach utilises goal-oriented error analysis using the concept of \emph{constitutive relation error} (CRE)~\cite{Pierre,CRE}. Originally,~\cite{Ladeveze} applied CRE to the evaluation of discretization errors in a finite element context, which is similar in terms of practicality and implementation to an equilibrated residual approach~\cite{Oden,Stein}. In a related approach, for linear elliptic diffusion problems,~\cite{Paladim} presents an approach to bound the error caused by the removal of stochastic inner features (homogenisation) without any integration of the fine‐scale features. The construction of the CRE method is based upon the satisfaction of the constitutive equation that guarantees the error in the energy norm bound for a linear QoI and increases its sharpness to estimate a more accurate error in the QoI. The finite element solution must be kinematically admissible, and the error is found by minimizing the potential energy in the energy norm. The static admissibility condition is obtained for the flux of the solution of second order divergence PDEs by minimizing the energy norm of the complementary energy theorem~\cite{Ladeveze,Ladeveze1}. Both potential energy and complementary energy theorems are helping to split up the error in the constitutive relation into two different minimization error measures, separately allocated to the field and flux. Both, complementary and classical potential energy contribute to the computation of the error in the energy norm for the CRE. CRE has the advantage that it does not require us to know the solution of the original model, and only requires the computation of the energy norms, which is fast, robust and accurate.

The idea of goal-oriented error estimation is to link the operator of the PDE with its adjoint operator, giving rise to a dual boundary value problem~\cite{Sueli}. The error between the recovered flux and finite element flux (stress in elasticity or electric displacement in electrostatics) for the simplified model in the energy norm bounds the difference between the QoI of the original and simplified problem. We assume that the boundary value problem has a unique solution, and that the PDE operator is self-adjoint such as the Laplacian or the Laplace-Beltrami operator.

In this paper we derive computationally affordable, strict upper and lower bounds (up to numerical approximation errors) for the simplification error for a QoI which is linear in the solution, based on goal-oriented error analysis and CRE, employing static and kinematic admissibility conditions. Errors resulting from linearizing a QoI are assumed to be negligible compared to the simplification error. The approach works with all feature types (internal, negative and positive), subject to different boundary conditions. We are able to study the simplification error in a localized subdomain and bound the error by only computing the norm in that subdomain. Our technique is demonstrated for electrostatics, but can be adapted to other common second order PDEs. The static admissibility condition needs to be satisfied for the flux of the finite element model with respect to the constitutive relation and the kinematic admissibility condition enforces the field (the solution of the PDE) to be equal to the Dirichlet value on the related boundary.

\section{Electrostatics}

We briefly review the governing equations for electrostatics. Assuming that we only deal with homogeneous, isotropic, non-dispersive, linear materials, Maxwell's equations simplify to Gauss' law~\cite{Jackson},
\begin{align}
 \n \cdot D & = \rho,\\
 D & = \varepsilon_r E, \label{eq:DE}
\end{align}
where $D$ is the electric displacement field, $\rho$ the charge density, $E$ the electric field and $\varepsilon_r$ the relative dielectric material constant, i.e.\ the permittivity expressed as a ratio relative to the permittivity of vacuum $\epsilon_0$. (We limit the discussion here to real-valued relative permittivity; in general in anisotropic materials the permittivity is a rank two tensor and is frequency dependent).  Introducing the electrostatic potential $\Phi$ via $E = -\nabla \Phi$, gives a \emph{Poisson equation:}
\begin{equation}
 \n^2 \Phi = - \rho/\varepsilon_r.\label{eq:PSE}
\end{equation}

This simplifies to a \emph{Laplace equation} in regions of space where $\rho \equiv 0$. For problems involving only insulators and conductors, $\rho$ is only non-zero on the boundary of the conductors (inside a conductor the voltage is constant, so $E$, $D$, and $\rho$ vanish, and there are no charges in an insulator). To turn Eq.~\eqref{eq:PSE} into a uniquely solvable boundary value problem, Dirichlet and/or Neumann boundary conditions are required. Then \emph{Green's identities} ensure that specifying a potential on a closed boundary (e.g.\ on conductors with different potentials), corresponding to a Dirichlet boundary condition defines a unique potential. Similarly, we obtain a unique potential by specifying the normal derivative of the potential (the electric field) on the boundary (corresponding to a given charge density) for Neumann boundary conditions. Thus, overall, this results in the boundary value problem for electrostatics:
\begin{equation}
\begin{aligned}
\n \cdot (\varepsilon_r \n \Phi) & = -\rho &&\text{ in }\Omega,\\
\Phi & = \Phi_D &&\text{ on }\Gamma_D,\\
\mathbf{n} \cdot (\varepsilon_r \n \Phi) & = \rho_N &&\text{ on }\Gamma_N,
\end{aligned}\label{eq:ge}
\end{equation}
where $\Omega$, $\Gamma_D$, $\Gamma_N$ are respectively the whole domain, and the Dirichlet and Neumann boundaries.

The weak form of the governing equation is derived by multiplying Eq.~\eqref{eq:ge} by an arbitrary test function $\Psi$ and integrating over $\Omega$:
\begin{equation}
\int_\Omega \varepsilon_r \n \Phi \cdot \n \Psi \;\mathrm{d}\Omega = \int_{\Gamma_N} (\mathbf{n} \cdot (\varepsilon_r \n \Phi)) \Psi \;\mathrm{d}\Gamma_N + \int_\Omega \rho \Psi \;\mathrm{d}\Omega. \label{eq:wf}
\end{equation}
As the potential $\Phi$ is smooth, $\Phi$ and $\Psi$ are in the Sobolev space $W^{1,p}$, which for $p=2$  is the Hilbert space $H^1(\Omega):=\{ \Phi \in L_2(\Omega) : \n \Phi \in L_2(\Omega) \}$ with the inner product $\langle \Phi,\Psi \rangle_{L_2} + \langle \n \Phi,\n \Psi \rangle_{L_2}$. Specific conditions (for capacitor simulations), prescribing that the value of the normal vector of the electrostatic potential at the outer boundary of the domain should always be zero, compel the charge density or the electrostatic potential derivative to be zero on the Neumann boundary. In this case the weak form simplifies to
\begin{equation}
B(\Phi,\Psi) = \ell(\Psi), \text{ where }
B(\Phi,\Psi) = \int_\Omega \varepsilon_r \n \Phi \cdot \n \Psi \;\mathrm{d}\Omega, \quad
\ell(\Psi) = \int_\Omega \rho \Psi \;\mathrm{d}\Omega\label{eq:2.8}
\end{equation}
are the bilinear and linear form respectively. We assume $\varepsilon_r$ is piecewise constant, such that we can partition $\Omega$ into subdomains $\Omega_l$, $\bigcup_{l=1}^{N} \Omega_l = \Omega$, where $\varepsilon_r$ has a constant permittivity $\epsilon_l$ on each $\Omega_l$.

To solve the weak form with finite element analysis, the function space of test and trial functions is discretized using a suitable basis $\{N_l\}$ such that $\Phi \approx \Phi_h = \sum_l \Phi_l N_l$. The $N_l$ are called shape functions, determined by the discretization of the domain and a class of basis functions. The elements must meet the minimum requirements for finite element simulation: continuity, consistency and completeness. For electrostatics, typically, arbitrary-order polynomial functions are well suited. Substituting the function approximations ($\Phi_h,\Psi_h$) into the weak form gives a linear equation $K \Phi_h = l$ where
\begin{align}
K_{ij} = \int_\Omega \varepsilon_r \n N_i \cdot \n N_j \;\mathrm{d}\Omega,\\
l_i = \int_{\Gamma_N}   (\mathbf{n} \cdot (\varepsilon_r \n \Phi)) N_j \;\mathrm{d}\Gamma_N + \int_\Omega \rho N_j \;\mathrm{d}\Omega.
\end{align}
The shape, hp-order and type of refinement of each element may cause noticeable differences in the solution~\cite{Monk03}. Here, for our numerical results, we refine the finite element mesh sufficiently and use polynomials of sufficient degree such that the discretization error is negligible for the problems, as verified by testing the convergence of the numerical solutions. The solution is obtained for homogeneous and inhomogeneous boundary conditions. It is important to reconstruct the affine function space in the particular case of inhomogeneous boundary conditions to obtain the solution making use of lifting for the Dirichlet boundary value, $\Phi = \Phi_0 + \Phi_D$~\cite{Demkowicz,Zienkiewicz}. $\Phi_0$ denotes the finite element solution w.r.t.\ the lifting ansatz. This gives the modified linear form and with it a new right hand side, $B(\Phi_0, \Psi) = \ell(\Psi) - B(\Phi_D, \Psi)$. Note that the solution is kinematically admissible, which is an important condition in the derivation of our simplification error estimate. The field is kinematically admissible when it becomes equal to the Dirichlet value on the boundary condition.

\section{Model Simplification Error Estimation} \label{General}

In this section we derive our goal-oriented approach toward estimating the model simplification error. Initially we explain the different possibilities of feature location and defeaturing. In the rest of the section, we detail our a posteriori goal-oriented simplification error estimation method with a QoI linear in the solution of the electrostatics problem. To bound the simplification error, we employ the constitutive relation error (CRE) technique~\cite{Pierre,CRE,Ladeveze1} to construct the bound in the energy norm. These will later on be applied to estimate the simplification error for internal features, and negative and positive boundary features for electrostatics problems.

\subsection{Defeaturing and Model Simplification}

\begin{figure}[t]
 \centering
 \begin{tabular}{cc}
 \includegraphics[width=0.3\textwidth]{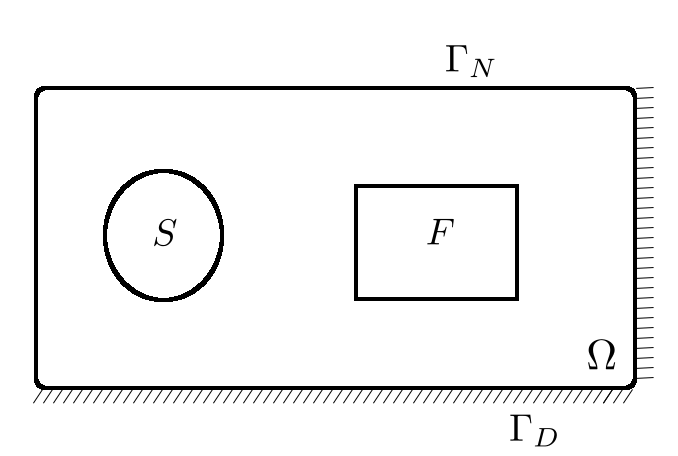}&
 \includegraphics[width=0.3\textwidth]{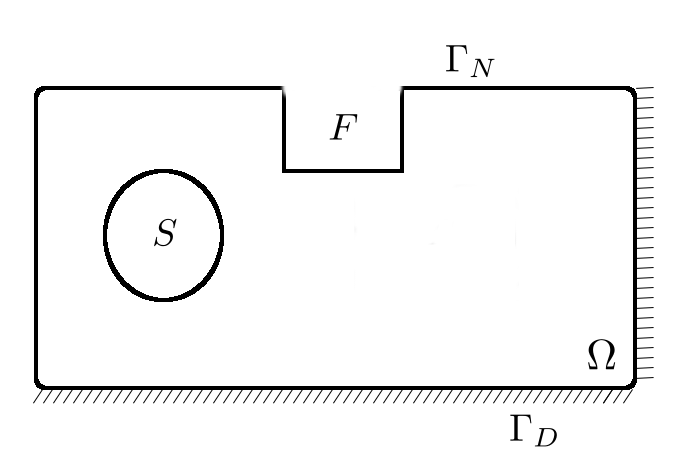}\\
 (a)&(b)\\
 \includegraphics[width=0.3\textwidth]{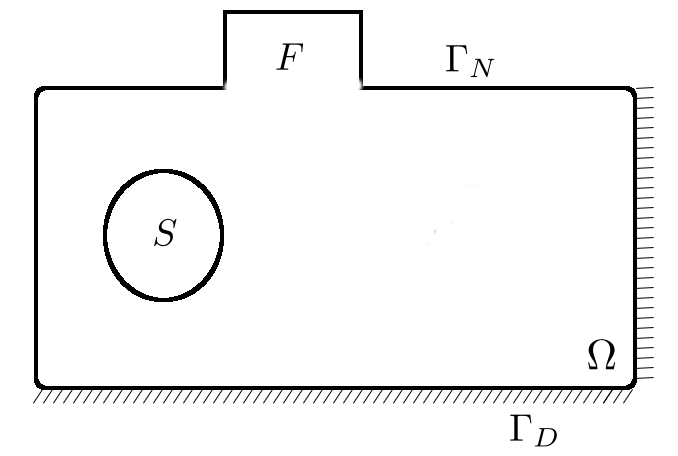}&
 \includegraphics[width=0.3\textwidth]{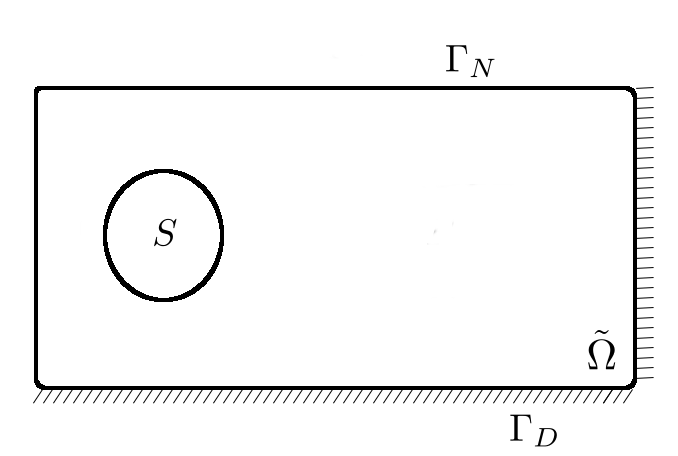}\\
 (c)&(d)
 \end{tabular}
 \caption{Models with (a) internal feature, (b) negative and (c) positive boundary feature and their simplified model (d) where $F$ is the feature domain, $\Gamma_D$ is part of boundary with Dirichlet boundary conditions, $\Gamma_N$ is part of boundary with Neumann boundary conditions, $S$ is the domain of interest for the QoI, $\Omega$ the domain of the original model and $\tilde{\Omega}$ the domain of the simplified model.}
 \label{fig:features}
\end{figure}

A feature is defined in general as a subset of a CAD model or even it could be a missing part related to the CAD model, typically associated with some semantic context. For example, holes, pockets, slots, etc. are referred to as features. For the purpose of this article, we do not consider specific feature types, but only refer to them as subsets of the boundary value problem domain $\Omega$ or a suitable extension of it within its span (for negative boundary features). The dimension $d$ (of the span) of $F$ is assumed to be equal to that of the domain $\Omega$. This leaves us with three general types:
\begin{itemize}
\item An \emph{internal feature} $F$, lies inside the domain, i.e. $F\subset\Omega$ with $F\cap\partial\Omega = \emptyset$. Its simplification consists of a change of the functions (e.g. material properties) defined on $F$ without changing the domain $\Omega$. For example, the feature can be a domain with material properties different from those of the surrounding domain. By simplifying the internal feature, the domain becomes more uniform as there are no discontinuities; see Figs.~\ref{fig:features}(a),(d).
\item A \emph{negative boundary feature} $F$ is an intrusion into the boundary of the geometry. In the original problem the negative feature can take the formation of a void without material, and it changes the boundary conditions by removing it. This means there is an extended domain $\tilde\Omega = \Omega \cup F$ in $\R^d$ such that $F\cap\partial\Omega \neq \emptyset$ and $F \cap \Omega^\circ = \emptyset$ (where $\Omega^\circ$ is the interior of $\Omega$); see Figs.~\ref{fig:features}(b),(d).
\item A \emph{positive boundary feature} is a protrusion on the boundary of the domain $\Omega$. It can contain a different material property to the surrounding domain, too. If it is removed, the boundary conditions will change subsequently next to removing whatever material property inside the feature. This means the domain $\Omega$ is partitioned into a feature domain $F$ and a remaining defeatured domain $\tilde\Omega$ such that $\Omega = \tilde\Omega \cup F$ where $\tilde\Omega$ and $F$ are connected by a common boundary section; see Figs.~\ref{fig:features}(c),(d).
\end{itemize}
Positive and negative features generally make the boundary more complex, meaning meshing is more expensive. The internal feature can also disturb the uniform mesh inside $\Omega$. To simplify the boundary, the feature can be either filled with the same material as the rest of domain if it is a negative feature or cut out of the boundary if it is a positive feature. An internal feature can be simplified by changing its material properties to match those of the surrounding domain. The internal feature can also be adjacent to the boundary if simplifying it does not affect the boundary of $\Omega$. In the following we discuss simplifying a single feature of one particular type.

\subsection{A Posteriori Error Estimation Approach}

Our approach is based on a posteriori goal-oriented error estimation, which was devised for estimating and subsequently reducing the simulation error for a particular QoI~\cite{Ming1,Oden-Prud,Oden-Prud2}. Goal-oriented error estimation provides a framework for relating the residual error (the main source of computational error) to the estimated QoI value. When considering simplification, the residual error becomes the difference between the finite element solutions for the original and simplified problems. To bound the error measures in the energy norm, we utilize a strategy based on CRE which leads to guaranteed bounds for the error in the QoI, providing a robust error estimator for the simplification error.

\subsubsection{Quantity of Interest (QoI) Error Estimation} \label{QoI}

Let $\Delta$ be a second order linear partial differential operator and $f$ be a sufficiently smooth function over $\Omega \subset \R^d$. We wish to solve the equation $\Delta u = f$ in $\Omega$ with homogeneous or inhomogeneous boundary conditions; we call this the original problem. Similarly let $\Delta\tilde{u} = \tilde{f}$ be the simplified version of the problem over another domain $\tilde\Omega$. Note that $\tilde{\Omega}$ may consist of sub-domains with different material constants from those in $\Omega$ and may have additional or removed domains depending on the defeaturing operation. In order for the differences between various functionals and the integrals employed by our estimation to be well defined, there is a need to define a \emph{simplified domain} $\widehat{\Omega}$. This domain must be compatible with $\Omega$. It is based on a feature domain $F$ and represents a modification of the defeatured domain $\tilde{\Omega}$ in order to make the computation of the bound possible. For the three different feature types we get the following cases:
\begin{enumerate}
\item[(i)] $\widehat{\Omega} = \tilde\Omega$ for an \textbf{internal feature}: the defeatured domain $\tilde{\Omega}$ and the simplified domain $\widehat{\Omega}$ are the same (but the functionals defined on them differ, e.g., due to different material properties).
\item[(ii)] $\widehat{\Omega} = \tilde\Omega \setminus F$ for a \textbf{negative boundary feature}: the simplified domain $\widehat{\Omega}$ is constructed from the defeatured domain $\tilde{\Omega}$ by removing the negative feature domain to enable the evaluation of the energy norm integrals.
\item[(iii)] $\widehat{\Omega} = \tilde\Omega$ for a \textbf{positive boundary feature}: this follows the same idea in the construction of $\widehat{\Omega}$ with the difference of adding the feature domain $F$, again to make the integrals over differences well defined.
\end{enumerate}
In the following we use $\Phi$ and related to refer to the fields, functionals and sets relating to the original problem and its domain $\Omega$; $\tilde\Phi$ and related to refer to the corresponding entities of the defeatured problem and its domain $\tilde\Omega$, and $\widehat{\Phi}$ and related to refer to the corresponding entities of the simplified problem and its domain $\widehat{\Omega}$. $\tilde{\Phi}$ and $\widehat{\Phi}$ both represent the simplified problem, but $\widehat{\Phi}$ is defined on the same domain than $\Phi$ to enable the comparison for the error bound calculation. The construction of $\widehat{\Omega}$ also includes the addition of a new boundary condition between $F$ and $\tilde{\Omega}$ that facilitates $\widehat{\Omega}$ to become compatible to the original domain $\Omega$. Details of their construction will be explained later when we discuss the different feature types. For the general derivation of the error bounds we will always refer to $\widehat{\Omega}$, $\widehat{\Phi}$ and related here.

Let $\widehat{\Phi}$  be the finite element solution of the boundary value problem Eq.~\eqref{eq:ge} over the simplified domain $\widehat{\Omega}$ and $\Phi$ be the solution of the original problem over $\Omega$, for the a posteriori goal-oriented error estimation~\cite{Oden,Oden-Prud}. Let $Q$ be the QoI for the problem, which must be linear in the solution. It is often an integral involving the solution over a subset of $\Omega$, e.g. $Q(\Phi) = \frac{1}{\mid S \mid} \int_{S} \ell(\Phi)\; \mathrm{d}\Omega$ where $S \subset \Omega$ is the domain of interest, and $\ell$ is the linear functional for the original problem. We wish to approximate
\begin{equation}
Q(\Phi) - Q(\widehat{\Phi}) = Q(\Phi - \widehat{\Phi}) = Q(e)  \text{ for }\;e = \Phi - \widehat{\Phi}.
\label{eq:QOI}
\end{equation}
The linearity of the QoI plays a significant role in the derivation of the simplification error estimate as it enables to take the difference between the QoIs for the simplified and original models. Note that in case the desired $Q$ is not linear, it can be linearized in certain cases, introducing further errors, which are not considered in detail here. Moreover, defeaturing must not affect the domain of interest $S$, i.e. $F \cap S = \emptyset$. This is justified as any features overlapping with $S$ are very likely to have a large effect on the QoI and so should not be considered for defeaturing.

We only use the solution of the simplified primal and dual boundary value problems to construct an upper and lower bound for the error in the QoI as explained now. The general weak form of the electrostatic PDE, given in Eq.~\eqref{eq:wf}, yields the weak forms of the original and simplified electrostatic problems,
\begin{eqnarray}
\label{eq:PF}
  B(\Phi,\Psi)  = \ell(\Psi),\quad \Phi \in U, \quad \forall \Psi \in V,\\
\label{eq:PS}
  \widehat{B}(\widehat{\Phi},\Psi)  = \widehat{\ell}(\Psi),\quad \widehat{\Phi} \in U, \quad \forall \Psi \in V,
\end{eqnarray}
respectively, where $U$ and $V$ are suitable function spaces for the trial and test functions, $\ell: V \to \R$ and $\widehat{\ell}: V \to \R$ are linear forms, and bilinear forms $B: U \times V \to \R$ and $\widehat{B}: U \times V \to \R$ are symmetric and positive definite, which define an inner product on $U$ and $V$. The \emph{energy norm} is defined as
\begin{equation}
  \| \Phi \|_{\varepsilon_r}:= \sqrt{B(\Phi,\Phi)} = \sqrt{\int_\Omega \varepsilon_r \n \Phi \cdot \n \Psi\; \mathrm{d}\Omega}.
\end{equation}
The test functions $\Psi$ for the primal and dual model are taken from
\begin{equation}
 V = \{ \Psi \in H^1(\Omega) \mid \Psi = 0 \text{ on } \Gamma_D \} \label{eq:test}.
\end{equation}
For the dual model, $\Psi = 0$  on $\Gamma_N$ where $\Gamma_N$ is the boundary of $\Omega$ with Neumann boundary conditions and $\Gamma_D$ the boundary of $\Omega$ with Dirichlet boundary conditions such that $\Gamma_N \cap \Gamma_D = 0$, $\Gamma_D \cup \Gamma_N = \partial \Omega$. As $B$ and $\widehat{B}$ are defined on the same domain, we can define the \emph{residual}
\begin{equation}
  R(\Psi) = B(\Phi-\widehat{\Phi},\Psi) = B(e,\Psi), \quad \forall \Psi \in V. \label{eq:Res}
\end{equation}

The choice of QoI depends on the engineering problem to be studied, the physics and governing equation. Here we choose the electric energy stored in the domain of interest $S$. It is the sum of all potential work that can be done in this domain by the electric field,
\begin{equation}\label{eq:q}
  q(\Phi) = \int_{S} D \cdot E \; \mathrm{d}\Omega = \int_{S} \varepsilon_r \n^2\Phi \; \mathrm{d}\Omega.
\end{equation}
This is a quadratic function, so we must linearize $q(\Phi)$ by replacing it with the first term of its Taylor expansion. Due to the linearity requirement this approximation is necessary, and we assume the nonlinear component is negligible compared to the error caused by simplification. The primal solution $\Phi$ is replaced by $\Psi$ (the test function) due to the equality of the function spaces for test and trial functions and so the linearized QoI becomes
\begin{equation}
  q^{'}(\Psi) = \int_{S} \varepsilon_r \n\Phi \cdot \n\Psi \; \mathrm{d}\Omega.
\end{equation}
However, $\Phi$, the solution of the original primal model, is unknown. So we assume it can be substituted by the solution of the simplified model $\widehat{\Phi}$, giving the linearized QoI of the simplified model,
\begin{equation}
  Q(\Psi) = \int_{S} \varepsilon_r \n\widehat{\Phi} \cdot \n\Psi \; \mathrm{d}\Omega \approx q(\Phi). \label{eq:qoi}
\end{equation}
While for electrostatic fields these assumptions are often fulfilled, this of course reduces the tightness of the bound and the guarantees for the upper and lower bound only hold for the linear QoI. With this approximation we may now use Eq.~\eqref{eq:QOI} to estimate the simplification error from the finite element solution of the simplified problem (solving the primal and dual simplified problems only). The benefit of this is to estimate the simplification error in a particular subset of the original model without expensive solving of the original problem.

\subsubsection{Dual (adjoint) Model}

Goal-oriented a posteriori error estimation requires to define an adjoint problem, which seeks a generalized Green's function associated with the QoI. It enables to localize the error to the domain of interest associated with the QoI. As $\Delta$ is a self-adjoint operator, we introduce the original and simplified dual problems to the primal problems in Eqs.~\eqref{eq:PF} and~\eqref{eq:PS}:
\begin{eqnarray}
\label{eq:DF}  B(\Psi,\Phi^\ast)  = Q(\Psi),\quad \Phi^\ast \in U, \quad \forall \Psi \in V, \\
\label{eq:DS}  \widehat{B}(\Psi,\widehat{\Phi}^\ast)  = Q(\Psi),\quad \widehat{\Phi}^\ast \in U, \quad \forall \Psi \in V,
\end{eqnarray}
where the respective bilinear forms are the same as in the primal problems with suitable boundary conditions. The choice of linear forms on the right hand side is the QoI, $Q$. This is quite often the suitable choice for many goal-oriented a posteriori error estimation problems.

The role of the dual model in goal-oriented error estimation is to relate the error in the QoI to the source of the error via setting the right hand side to the QoI over the domain of interest $S$. The dual boundary value problem should always be homogeneous, so it has a dual solution, $\Phi^\ast$, for the non-zero right hand side for the governing equation. It is derived similarly to Eq.~\eqref{eq:ge}, except for the fact that the right hand side is given by the QoI and Dirichlet and Neumann boundary conditions are set to zero:
\begin{equation}
\begin{aligned}
D^\ast - \widehat{D}_s &= -\varepsilon_r \n\Phi^\ast &&\text{ in } S,\\
D^\ast &= -\varepsilon_r \n\Phi^\ast &&\text{ in } \Omega/ S, \\
\n \cdot (- \varepsilon_r \n\Phi^\ast) &= Q(\Psi) &&\text{ in } \Omega, \\
\Phi^\ast &= 0 &&\text{ on } \Gamma_D,\\
\mathbf{n} \cdot ( \varepsilon_r \n\Phi^\ast) &= 0 &&\text{ on } \Gamma_N,
\end{aligned}\label{eq:dual}
\end{equation}
where $D^\ast$ is the dual electric displacement and $\widehat{D}_s$ is the flux of the solution of the simplified primal model in the region of interest $S$. Like in the primal case, the original and simplified models give dual solutions $\Phi^\ast$ and $\widehat{\Phi}^\ast$ respectively. For the dual problem, the Neumann and Dirichlet boundary conditions are equal to zero while the right hand side of the PDE is no longer zero. Note that the solution of the dual problem satisfies the kinematic admissibility condition on the Dirichlet boundary due to the homogeneous boundary condition. But it is difficult to target the equilibrium constitutive equation to satisfy the static admissibility condition for the flux.

The bilinear form derived from the original dual problem for Eq.~\eqref{eq:dual} is
\begin{equation}
B(\Phi^\ast , \Psi) = \int_\Omega \varepsilon_r \n\Phi^\ast \cdot \n\Psi \;\mathrm{d}\Omega = -\int_{S} \widehat{D}_s \cdot \n\Psi \; \mathrm{d}\Omega = \int_{S} \varepsilon_r \n\widehat{\Phi}^\ast \cdot \n\Psi \; \mathrm{d}\Omega.
\end{equation}
The function space $V$ of the test function $\Psi$ is the same than for the primal model, $V = \{ \Psi \in H^1 \mid \Psi = 0 \text{ on } \Gamma_D \text{ and } \Gamma_N \}$. Similarly, the bilinear form of the weak form for the simplified dual problem is
\begin{equation}
B(\widehat{\Phi}^\ast,\Psi) = \int_{\widehat{\Omega}} \varepsilon_r \n \widehat{\Phi}^\ast \cdot \n\Psi \; \mathrm{d}\Omega = \int_{S} \varepsilon_r \n\widehat{\Phi}^\ast \cdot \n\Psi \; \mathrm{d}\Omega.
\end{equation}

As we have chosen the same functional spaces for primal and dual models and $e = \Phi - \widehat{\Phi}$ and $e^\ast = \Phi^\ast - \widehat{\Phi}^\ast$ can be considered particular test functions, we get
\begin{equation}
  Q(e) = B(e,\Phi^\ast-\widehat{\Phi}^\ast) + B(e,\widehat{\Phi}^\ast) = B(e,e^\ast) + R(\widehat{\Phi}^\ast).
\end{equation}
Here $R(\widehat{\Phi}^\ast)$ is the residual for the simplified dual problem, similar to Eq.~\eqref{eq:Res}. Hence,
\begin{equation}
  Q(e) - R(\widehat{\Phi}^\ast) = B(e,e^\ast).
\end{equation}
The Cauchy-Schwarz inequality then gives the bound,
\begin{equation} \label{eq:berror}
  |Q(e) - R(\widehat{\Phi}^\ast)| \leq \sqrt{B(e^\ast,e^\ast)} \sqrt{B(e,e)}
  = \| \bigtriangledown \mathit{e}^\ast \|_{\varepsilon_r} \| \bigtriangledown \mathit{e} \|_{\varepsilon_r} \;
  \leq \nu^\ast \nu,
\end{equation}
where $\|.\|_{\varepsilon_r}$ denotes the energy norm of the error over the domain $\Omega$, and $\nu$ and $\nu^\ast$ are global estimates for the norms of the error of primal and dual simplified solutions. The bounds in Eq.~\eqref{eq:berror} are not computable as the energy norm of the exact error fields is not available. Instead, we calculate bounds, $\nu$ and $\nu^\ast$, for these quantities. This is elaborated in Section~\ref{CRE} where the CRE is employed to bound the energy norms. Provided that the bounds are sufficiently sharp and can be computed with reasonable effort, they bound the simplification error in the energy norm from the finite element solution of the primal and dual simplified model. Hence, the construction of the bound components must always satisfy the principle of virtual work and the constitutive equation. Thus, the admissibility conditions must be satisfied.

\subsubsection{Constitutive Relation Error}\label{CRE}

We take advantage of the \emph{constitutive relation error} (CRE) to bound the the energy norms $\nu$ and $\nu^\ast$ for primal and dual models. It only requires to employ the admissibility conditions. CRE provides a bound that is conceptually simple to understand, implement and control. It constructs a recovered electrostatic displacement that is statically admissible, or equilibrated. CRE applies the kinematic admissibility condition on the field which must be verified for the finite element model and its boundary conditions. The distance calculated in the energy norm between the recovered flux (electrostatic displacement) and simplified finite element electrostatic displacement is a bound for the simplification error.

The solution $\widehat{\Phi}$ of the simplified model is kinematically admissible, in other words the field $\widehat{\Phi}$ meets the Dirichlet boundary conditions. It will be sought in
\begin{equation}
  U = \{ \widehat{\Phi} \in H^1(\Omega) \; | \; \widehat{\Phi} = \Phi_D \text{ on } \Gamma_D \}.
\end{equation}
The flux, $\widehat{D}$, must be statically admissible, i.e.
\begin{equation}
 \int_{\widehat{\Omega}} \widehat{D} \cdot \n\Psi \; \mathrm{d}\Omega = \int_{\widehat{\Omega}} \rho \Psi \; \mathrm{d}\Omega + \int_{\Gamma_N} \mathbf{n} \cdot (\varepsilon_r \n\Phi) \Psi \; \mathrm{d}\Gamma, \label{eq:stat-con}
\end{equation}
which means that the constitutive equation must be always satisfied.

The relative permittivity $\varepsilon_r$ in Eq.~\eqref{eq:ge} is assumed to be a piecewise constant and scalar, naturally giving rise to a partition of the domain into domains of constant permittivity. For the original model, this means in the simplest case,
\begin{equation}
 \varepsilon_r(x) =
  \begin{cases}
    \epsilon_R  & \text{for } x \in \Omega\setminus F,\\
    \epsilon_F  & \text{for } x \in F
  \end{cases}
\end{equation}
where $\epsilon_R,\epsilon_F \in \R^+_0$. We derive the error bounds for these constants, enabling us to construct the CRE error bounds for all concerned feature types. This is used later in the computation and evaluation of the error bounds for different feature types.

To make Eq.~\eqref{eq:berror} practically useful, we must find the bounds $\nu$ and $\nu^\ast$ for $B(e,e)$
and $B(e^\ast,e^\ast)$ respectively. For $B(\Phi-\widehat{\Phi},\Phi-\widehat{\Phi}) = B(e, e) \leq \nu^2$, noting that $\widehat{D} = -\widehat{\varepsilon}_R \n \widehat{\Phi}$, let
\begin{equation}
\begin{aligned}
  \nu^2 &= \| \widehat{D} + \epsilon_R \n\widehat{\Phi} \|^2_{\epsilon_R^{-1}} + \| \widehat{D} + \epsilon_F \n\widehat{\Phi} \|^2_{\epsilon_F^{-1}} \label{eq:nunor} \\
        &= \int_{\widehat{\Omega}/F} (\widehat{D} + \epsilon_R \n\widehat{\Phi}) \epsilon_R^{-1} (\widehat{D} + \epsilon_R \n\widehat{\Phi}) \; \mathrm{d}\Omega\\
        &\qquad {} + \int_{F} (\widehat{D} + \epsilon_F \n\widehat{\Phi}) \epsilon_F^{-1} (\widehat{D} + \epsilon_F \n\widehat{\Phi}) \; \mathrm{d}\Omega.
\end{aligned}
\end{equation}
This definition requires to calculate two norms, one over the feature domain $F$ with relative permittivity $\epsilon_F$, and the other one covers the error measure in the energy norm for the rest, $\Omega \setminus F$, with relative permittivity $\epsilon_R$. Note that in order to eliminate additional numerical approximations and return the exact value of the integration in the feature domain, we keep the mesh in the feature domain unchanged from the original to the simplified model. We choose to compute the energy norms in terms of electric displacement, which is a linear function related to the electrostatic potential. The electric displacement of the simplified problem, $\widehat{D}$, follows the same formula as the electric displacement $D$ defined in Eq.~\eqref{eq:DE} with a different finite element solution for the simplified model in the domain $\widehat{\Omega}$. $\widehat{D}$ helps to construct the flux to distinguish the feature domain $F$ from other parts of the original domain $\Omega$ for the computation of the error bounds.

By subtracting $0 = D + \varepsilon_r \n\Phi$ in the definition of $\nu^2$ we get
\begin{equation}
\begin{aligned}
  \nu^2 & = \| \widehat{D} - D - \epsilon_R \n\Phi + \epsilon_R \n\widehat{\Phi} \|^2_{\epsilon_R^{-1}} + \| \widehat{D} - D - \epsilon_F\n\Phi + \epsilon_F \n\widehat{\Phi} \|^2_{\epsilon_F^{-1}}\\
        &= \| (\widehat{D} - D) + \epsilon_R (\n\widehat{\Phi} - \n\Phi) \|^2_{\epsilon_R^{-1}} + \| (\widehat{D} - D) + \epsilon_F (\n\widehat{\Phi} - \n\Phi) \|^2_{\epsilon_F^{-1}}\\
        &= \int_{\widehat{\Omega}/F} (\widehat{D} - D) \epsilon_R^{-1} (\widehat{D} - D) \; \mathrm{d}\Omega + \int_F (\widehat{D} - D) \epsilon_F^{-1} (\widehat{D} - D) \; \mathrm{d}\Omega \\
        &\qquad {} + \int_{\widehat{\Omega}/F} (\n\widehat{\Phi} - \n\Phi) \epsilon_R^{-1} (\n\widehat{\Phi} - \n\Phi) \; \mathrm{d}\Omega\\
        &\qquad {} + \int_F (\n\widehat{\Phi} - \n\Phi) \epsilon_F^{-1} (\n\widehat{\Phi} - \n\Phi) \; \mathrm{d}\Omega  + 2 \underbrace{\int_{\widehat{\Omega}} (\widehat{D} - D)(\n\widehat{\Phi} - \n\Phi) \;\mathrm{d}\Omega}_{=0}.
\end{aligned} \label{eq:nu}
\end{equation}
Eq.~\eqref{eq:ge} follows from Eq.~\eqref{eq:DE} by taking the gradient under the divergence free condition. In combination with the weak form in Eq.~\eqref{eq:wf}, we get for $D$ and $\widehat{D}$ respectively:
\begin{align}
  \int_\Omega D \cdot \n\Psi \; \mathrm{d}\Omega &= \ell(\Psi), \label{eq:D} \\
  \int_{\widehat{\Omega}} \widehat{D} \cdot \n\Psi \; \mathrm{d}\Omega &= \ell(\Psi), \label{eq:hatD}
\end{align}
where $\Psi$ is substituted for the test function in the bilinear and linear forms. Subtracting Eq.~\eqref{eq:D} from Eq.~\eqref{eq:hatD} gives
\begin{equation}
  \int_{\widehat{\Omega}} (\widehat{D} - D) \cdot \n\Psi \; \mathrm{d}\Omega = \ell(\Psi) - \ell(\Psi) = 0. \label{eq:sub}
\end{equation}
This means the divergence free electric displacement $D$ and $\widehat{D}$ satisfy the static admissibility condition. $\widehat{\Phi}$ is kinematically admissible because it is equal to $\Phi$ on $\Gamma_D$. The divergence free condition for $\mathit{D}$ gives Eq.~\eqref{eq:sub}, leading to an orthogonality condition between statically and kinematically admissible variables by setting $\Psi = \Phi -\widehat{\Phi}$:
\begin{equation}
 \int_{\widehat{\Omega}} (D - \widehat{D})\cdot (\n\Phi - \n\widehat{\Phi}) \; \mathrm{d}\Omega = 0 \; \text{ as } \Phi - \widehat{\Phi} = 0 \text{ on } \Gamma_D.
\end{equation}
So the last term of the expanded Eq.~\eqref{eq:nu} is indeed zero and Eq.~\eqref{eq:nu} yields
\begin{equation}
  \nu^2 = \underbrace{\| \widehat{D} - D \|^2_{\epsilon_R^{-1}} + \| \widehat{D} - D \|^2_{\epsilon_F^{-1}}}_{\geq 0} + \| \n \widehat{\Phi} - \n\Phi \|^2_{\varepsilon_r} \label{eq:distance},
\end{equation}
which implies
\begin{equation}\label{eq:nub}
  \|e\|^2_{\varepsilon_r} = B(e,e) = \|  \n \Phi - \n \widehat{\Phi}\|^2_{\varepsilon_r} \leq \nu^2.
\end{equation}
Note that the closeness of the bound depends on the value of $\|\widehat{D}-D\|^2_{\epsilon_R^{-1}} + \|\widehat{D}-D\|^2_{\epsilon_F^{-1}}$. If the simplification only has a minor effect on the finite element solution, $\nu$ is nearly equal to the energy norm and hence provides a good estimate. If the effect is larger, we are overestimating it stronger.

Similarly, the CRE $\nu^\ast$ is derived from the simplified dual model solution $\widehat{\Phi}^\ast$:
\begin{equation}\label{eq:erdua}
\begin{aligned}
 (\nu^\ast)^2 &= \| \widehat{D}^\ast + \epsilon_R \n\widehat{\Phi}^\ast \|^2_{\epsilon_R^{-1}} + \| \widehat{D}^\ast + \epsilon_F \n\widehat{\Phi}^\ast \|^2_{\epsilon_F^{-1}}\\
              &= \int_{\widehat{\Omega}/F} (\widehat{D}^\ast + \epsilon_R \n\widehat{\Phi}^\ast) \epsilon_R^{-1} (\widehat{D}^\ast + \epsilon_R \n\widehat{\Phi}^\ast) \; \mathrm{d}\Omega\\
              &\qquad {} + \int_F (\widehat{D}^\ast + \epsilon_F \n\widehat{\Phi}^\ast) \epsilon_F^{-1} (\widehat{D}^\ast + \epsilon_F \n\widehat{\Phi}^\ast) \; \mathrm{d}\Omega,
\end{aligned}
\end{equation}
where $\widehat{D}^\ast$  is the flux of the solution of the simplified dual model. As before, similar to Eq.~\eqref{eq:nub}, the energy norm of the error, this time of the dual problem, is bounded,
\begin{equation}
  \| e^\ast \|^2_{\varepsilon_r} = B(e^\ast,e^\ast) = \| \n \Phi^\ast - \n \widehat{\Phi}^\ast\|^2_{\varepsilon_r} \leq (\nu^\ast)^2.
\end{equation}

After the computation of the distance error norms for primal and dual models, $\nu$ and $\nu^\ast$, we link the CRE error norm constructions to the QoI by the computation of the residual. Eq.~\eqref{eq:Res} gives
\begin{equation}\label{eq:residual}
\begin{aligned}
R(\widehat{\Phi}^\ast) & = B(e,\widehat{\Phi}^\ast) = B(\Phi,\widehat{\Phi}^\ast) - B(\widehat{\Phi},\widehat{\Phi}^\ast) = \ell(\widehat{\Phi}^\ast) - B(\widehat{\Phi},\widehat{\Phi}^\ast)\\
        & = \int_{\Gamma_N} (\mathbf{n} \cdot (\varepsilon_r \n\Phi))\widehat{\Phi}^\ast \; \mathrm{d}\Gamma + \int_\Omega \rho \widehat{\Phi}^\ast \; \mathrm{d}\Omega - \int_\Omega \varepsilon_r \n\widehat{\Phi} \cdot \n\widehat{\Phi}^\ast \; \mathrm{d}\Omega.
\end{aligned}
\end{equation}
This is the residual error of original versus simplified model. This residual equation is valid if and only if the bilinear form $B$ is self-adjoint, positive-definite and symmetric.

Finally, $R(\widehat{\Phi}^\ast)$, $\nu$ and $\nu^\ast$ in combination give the general form of the upper and lower bounds for the simplification error,
\begin{align}
 R(\widehat{\Phi}^\ast) - \nu \nu^\ast \leq &Q(e) \leq R(\widehat{\Phi}^\ast) + \nu \nu^\ast,\\
 Q(\widehat{\Phi}) + R(\widehat{\Phi}^\ast) - \nu \nu^\ast \leq &Q(\Phi) \leq Q(\widehat{\Phi}) + R(\widehat{\Phi}^\ast) + \nu \nu^\ast. \label{eq:gb}
\end{align}
In the following sections we apply these error bounds to various numerical problems for different feature types. The error bounds must be adapted to internal, negative and positive features for electrostatics problems and we show their performance in practical settings.

\section{Internal Features} \label{Inner}

In this section we apply the simplification error analysis to internal features in various capacitor problems. We first show how to adapt the theory and bounds to this specific problem and then present numerical results.

\subsection{Internal Feature Simplification Error Estimation}

An internal feature is a region $F \subset \Omega$ with $F \cap \d\Omega = \emptyset$ where the simplification consists of a change of the material properties in $F$, usually setting it to the material properties of the surrounding domain. This means the boundary conditions of the original versus simplified model remain intact and $\Omega = \tilde\Omega = \widehat{\Omega}$. The governing equations for the finite element simulation are given by Eq.~\eqref{eq:ge}.

\begin{figure}[t]
 \centering
 \includegraphics[width=.3\textwidth]{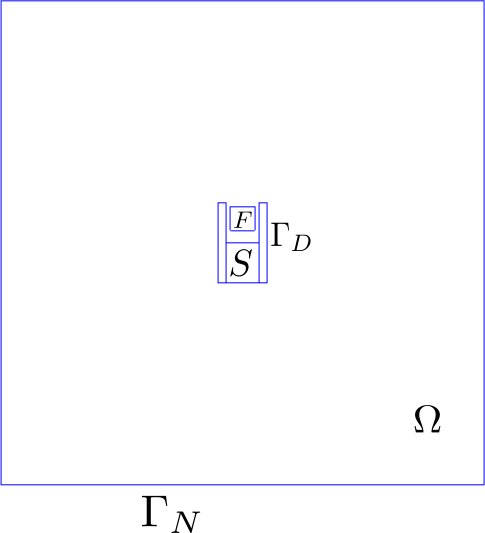}
 \caption{Parallel plate capacitor model for electrostatics simulation with domain of interest $S$ and feature $F$ in domain $\Omega$ with Dirichlet boundary $\Gamma_D$ and Neumann boundary $\Gamma_N$.}
 \label{fig:capacitor-feature}
\end{figure}

We consider the example of a parallel plate capacitor in Fig.~\ref{fig:capacitor-feature} where the feature $F$ and the domain of interest $S$ lie between the two capacitor plates with Dirichlet boundary conditions on $\Gamma_D$ and Neumann boundary conditions on the far boundary $\Gamma_N$ surrounding the capacitor. Note that, as the capacitor plates are conductors, the potential is constant (but different) on each plate boundary and we set it to $+220V$ and $-220V$ respectively. This produces a symmetric distribution of the field inside the capacitor. For the Neumann boundary, the far boundary has zero charge density, $\rho|_{\Gamma_N} = 0$, and the normal component of the electric displacement on the boundary is zero, $\mathbf{n} \cdot (\varepsilon_r \n \Phi) = 0$. Using our linear QoI approximating the electrostatic energy in the capacitor can be utilized to show the effect that a dielectric material has on the capacitance. We provide an upper and lower bound for the difference between the QoI of the original and simplified model.

Eq.~\eqref{eq:gb} with Eqs.~\eqref{eq:nunor},~\eqref{eq:erdua} and~\eqref{eq:residual} provide the bound for the simplification error of our QoI given by Eq.~\eqref{eq:qoi}. The residual Eq.~\eqref{eq:residual} must be adapted to the specific linear form $\ell(\Psi)$ of the governing Eq.~\eqref{eq:ge}. In this problem, there is no charge density in the domain $\tilde{\Omega}$, $\rho = 0$, so $\ell(\Psi) = 0$ in Eq.~\eqref{eq:residual}. The relative permittivity in $F$ is constant, i.e. $\varepsilon_r(x) = \epsilon_F$ for $x \in F$, but different from the surrounding permittivity in $\Omega \setminus F$. For the simplified problem, $\tilde\varepsilon_r$ is constant on the whole domain $\tilde{\Omega}$. Because the Dirichlet boundary condition is not zero, it needs to be shifted to the right hand side via lifting the boundary condition. This results in the lifted weak form,
\begin{equation}
  B(\Phi_0,\Psi) = \int_\Omega \varepsilon_r \n \Phi_0 \n \Psi \;\mathrm{d}\Omega = -\int_{\Gamma_D} \varepsilon_r \n \Phi_D \n \Psi \;\mathrm{d}\Omega = - B(\Phi_D,\Psi),
   \label{eq:lslift}
\end{equation}
where $\Phi_0$ is the solution. The simplified problem given by $\tilde{B}$, $\tilde{\ell}$ has a constant relative permittivity $\tilde\varepsilon_r \equiv \epsilon_R$  and $\Phi,\tilde{\Phi} \in H^1(\Omega)$ and $\Psi,\tilde{\Psi} \in H^1_0(\Omega)$.

To bound $Q(e)$, the solution $\tilde{\Phi}$ of the simplified problem is substituted into Eq.~\eqref{eq:nunor} to determine $\nu$ and for $\tilde{D}$ we get
\begin{equation}
  \int_{\tilde{\Omega}} (\n \cdot \tilde{D}) \Psi \; \mathrm{d}\Omega = 0.
\end{equation}
By noting that the Neumann boundary condition is zero and the test function $\Psi$ vanishes on $\Gamma_D$, integrating the above by parts gives
\begin{equation}
 \underbrace{\int_{\tilde{\Omega}} \n (\tilde{D}\Psi) \; \mathrm{d}\Omega}_{=\int_{\partial{\tilde{\Omega}}} (\mathbf{n} \cdot \tilde{D})\Psi \; \mathrm{d}\Gamma = 0} -
 \int_{\tilde{\Omega}} \tilde{D} \n\Psi \; \mathrm{d}\Omega = 0.\label{eq:adm}
\end{equation}
Using Eq.~\eqref{eq:sub}, subtracting the fluxes of the original and simplified models, means the static admissibility condition is fulfilled:
\begin{align}
 \int_{\widehat{\Omega}} (D-\tilde{D})\cdot\nabla\Psi \; \mathrm{d}\Omega = 0 \quad \forall \Psi \in H^1_0.
\end{align}
Eq.~\eqref{eq:distance} now lets us compute $\nu$ from the solution of the simplified model:
\begin{equation}
\begin{aligned}
\nu^2 =& \parallel \tilde{P} + \epsilon_R \nabla\tilde{\Phi} \parallel^2_{{\epsilon_R}^{-1}} + \parallel \tilde{D} + \epsilon_F \nabla\tilde{\Phi} \parallel^2_{{\epsilon_F}^{-1}}\\
  =& \underbrace{\int_{\tilde{\Omega}/F} (\tilde{D} + \epsilon_R \n \tilde{\Phi}){\epsilon_R}^{-1}(\tilde{D} + \epsilon_R \n \tilde{\Phi}) \; \mathrm{d}\Omega}_{ = 0}
   + \int_{F} (\tilde{D} + \epsilon_F \n \tilde{\Phi}){\epsilon_F}^{-1}(\tilde{D} + \epsilon_F \n \tilde{\Phi}) \; \mathrm{d}\Omega.
\end{aligned}
\end{equation}
$\nu$ depends on the integration over two separate domains, $F$ and $\Omega\setminus F$. The term for $\Omega \setminus F$ is zero due to the construction of $\tilde{D}$: field and flux only differ by the material property of the domain that they belong to. Technically, $\tilde{D}$ ensures that the solution of the simplified model, $\tilde{\Phi}$, and flux in $\Omega\setminus F$ are identical with the recovered flux and field in $\nu$. Therefore we only need to consider the feature domain $F$ with the different relative permittivities. Hence, from Eqs.~\eqref{eq:nu} and~\eqref{eq:distance},
\begin{equation}
\begin{aligned}
\nu^2 &=  \int_{F} (\epsilon_R \nabla\tilde{\Phi} - \epsilon_F \n \tilde{\Phi}){\epsilon_F}^{-1}(\epsilon_R \nabla\tilde{\Phi} - \epsilon_F \n \tilde{\Phi}) \; \mathrm{d}\Omega \\
 &= \int_{F} \frac{({\epsilon}_R - \epsilon_F)^2}{\epsilon_F} \n \tilde{\Phi} \n \tilde{\Phi} \; \mathrm{d}\Omega.
\end{aligned}\label{eq:prinormerr}
\end{equation}

Similarly, we calculate $\nu^\ast$ via Eq.~\eqref{eq:erdua} for the dual model, except that it now depends on the solution $\tilde{\Phi}^\ast$ of the dual simplified model:
\begin{equation}
\begin{aligned}
(\nu^\ast)^2 &= \parallel \tilde{D} + \epsilon_R \nabla\tilde{\Phi}^\ast \parallel^2_{{\epsilon_R}^{-1}} + \parallel \tilde{D} + \epsilon_F \nabla\tilde{\Phi}^\ast \parallel^2_{{\epsilon_F}^{-1}} \\
&= \underbrace{\int_{\widehat{\Omega}/F} (\widehat{D} + \epsilon_R \nabla\widehat{\Phi}^\ast) {\epsilon_R}^{-1} (\widehat{D} + \epsilon_R \nabla\widehat{\Phi}^\ast) \; \mathrm{d}\Omega}_{=0}\\
&\qquad {} + \int_F (\widehat{D} + \epsilon_F \nabla\widehat{\Phi}^\ast) {\epsilon_F}^{-1} (\widehat{D} + \epsilon_F \nabla\widehat{\Phi}^\ast) \; \mathrm{d}\Omega.
\end{aligned}
\end{equation}
The integration over $\widehat{\Omega}\setminus F$ is zero because the flux remains the same between the two models in this domain, so
\begin{equation}
  (\nu^\ast)^2 = \int_{F} \frac{({\epsilon}_R - \epsilon_F)^2}{\epsilon_F} \n \tilde{\Phi}^{\ast} \n \tilde{\Phi}^{\ast} \; \mathrm{d}\Omega.
\end{equation}

As for the capacitor, $\ell(\Psi) = 0$, the residual error $R(\tilde\Phi^\ast)$ given by Eq.~\eqref{eq:residual} simplifies to $R(\tilde\Phi^\ast) = -B(\tilde{\Phi},\tilde{\Phi}^{\ast})$. As $\tilde{\Phi}^{\ast}$ vanishes on the Dirichlet boundary, we have
\begin{multline}
B(\tilde{\Phi}^{\ast},\tilde{\Phi}) = \int_{F} \epsilon_F \n\tilde{\Phi}\n\tilde{\Phi}^{\ast} \mathrm{d}\Omega + \int_{\Omega/ F} \epsilon_R \n\tilde{\Phi}\n\tilde{\Phi}^{\ast} \mathrm{d}\Omega\\
= \int_{F} (\epsilon_F + \epsilon_R - \epsilon_R) \n\tilde{\Phi}\n\tilde{\Phi}^{\ast} \mathrm{d}\Omega + \int_{\Omega/ F} \epsilon_R \n\tilde{\Phi}\n\tilde{\Phi}^{\ast} \mathrm{d}\Omega\\
= \underbrace{\int_\Omega \epsilon_R \n \tilde{\Phi}\n\tilde{\Phi}^{\ast} \mathrm{d}\Omega}_{=0} + \int_{F} (\epsilon_F - \epsilon_R)\n\tilde{\Phi}\n\tilde{\Phi}^{\ast} \mathrm{d}\Omega.
\end{multline}
The first term in the last line above is zero because the bilinear form of the simplified model satisfies the general boundary value problem and its boundary conditions. Hence,
\begin{equation}\label{eq:finres}
R(\tilde\Phi^\ast) = \int_{F} (\epsilon_R - \epsilon_F)\n\tilde{\Phi}\n\tilde{\Phi}^{\ast} \mathrm{d}\Omega.
\end{equation}

With this we have all components to calculate the bounds for Eq.~\eqref{eq:gb}. The integrals are evaluated on the same mesh before and after defeaturing. This is to avoid additional numerical errors from approximating the integrals over the various domains by keeping the meshing in the feature domain consistent between original and simplified models.

\subsection{Numerical Results for Internal Features}

Here we present numerical results to evaluate the effectivity of the bounds using a simple capacitor model and simplifying different feature domains. \cite{Pierre, Ladeveze1}~have already proven that the CRE bounds are guaranteed and tight bounds can be achieved. We applied the CRE method in goal-oriented error estimation, which consequently makes the error bounds in the energy norm more accurate by localising the error computation to a sub-set of $\Omega$. The error estimation is implemented as C++ plugin for the finite element solver NGSsolve~\cite{Ngsolve} with the Netgen mesh generator~\cite{Netgen} to calculate and render the solutions. We show how well our method is able to estimate the simplification error. For this we study the quality of the bounds using an effectivity index comparing the bounds with the actual simplification error:
\begin{equation}
  \mathfrak{I} = \frac{|Q(\Phi)| + |U-L|}{|Q(\Phi)|} = 1 + \frac{|U-L|}{|Q(\Phi)|},
\end{equation}
where $U$ is our upper and $L$ our lower bound. The closer the effectivity index is to one, the better the bounds are.

In the following experiments, the geometry of the capacitor model and associated boundary conditions are described. We first show a practical case of a parallel plate capacitor with pyrex dielectric and a sodium contamination. The simplification error in terms of stored electrostatic energy in the domain of interest $S$ is then investigated for two different simplification techniques: (i) removal of a fixed feature domain with different dielectric material properties, different from the surrounding domain; (ii) removal of a feature that grows in size with fixed dielectric constant, different from the surrounding domain, is removed. In both cases the feature is removed by setting the dielectric constant inside its domain to the dielectric constant of the surrounding domain.

\subsubsection{Parallel Plate Capacitor with Glass Dielectric Material}

\begin{figure}[t]
 \centering
 \hfil\includegraphics[width=.4\textwidth]{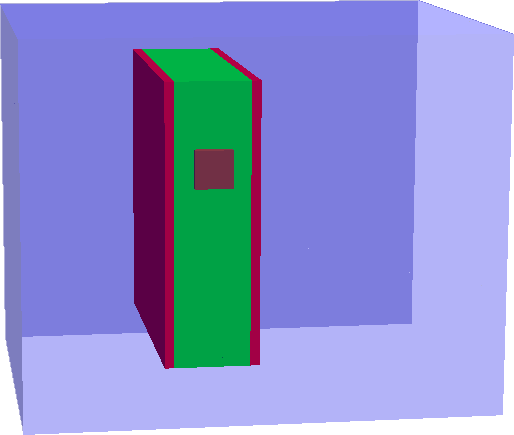}\hfil
 \includegraphics[width=.4\textwidth]{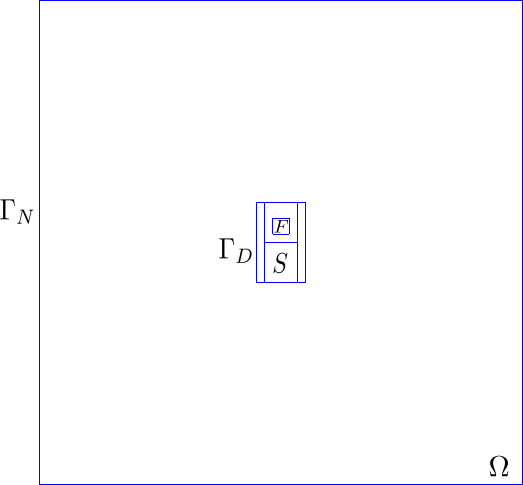}\hfil
 \caption{Left: 3D glass capacitor model with two metal conductor plates in red and. The whole box of the capacitor is in is filled with air (only partially shown to have a larger view of the capacitor itself). The dielectric material is in green (glass - Corning 7740) with a Sodium contamination in brown. Right: 2D model of the glass capacitor, where $\Omega$ is the whole domain. The Dirichlet boundary $\Gamma_D$ is at the conductor plates. The outer box of the capacitor is the Neumann boundary $\Gamma_N$. The domain of interest $S$ is the lower half between the capacitor plates. The sodium contamination is the domain $F$.}
 \label{fig:Capacitor-practical-3D}
\end{figure}

The first example model is a parallel plate capacitor shown in Fig.~\ref{fig:Capacitor-practical-3D} with the domain of interest $S$ and the feature domain $F$ being simplified lying between the two capacitor plates. There are many shapes for practical capacitors. However, they all consist of at least two electrical plates separated by a dielectric. The whole model domain for the capacitor is a square of $6$cm length. The conductor plates have a rectangular shape, each of the same width and length of $0.1$cm and $1$cm respectively, and they are placed in the centre of the domain. The conductor plates are $0.4$cm apart of each other. They are given a voltage of $+220$V and $-220$V. The right plate is the positive conductor. The capacitor is surrounded by air. The dielectric material between the capacitor plates stores the electric energy by becoming polarized, determining the capacitance. In order to maximize the capacitance, the dielectric material should have the highest possible permittivity. The dielectric material in our model is Corning 7740 (pyrex), a glass wafer. Glass provides reliable and stable performance and operates in a wide range of temperatures. The relative permittivity of pyrex is $4.6$~\cite{Nalwa}. When manufacturing this material, there is a probability of it being contaminated by sodium. Sodium contamination can be deleterious to the electrical properties of pyrex structures. The permittivity of sodium, higher than pyrex, can cause damage to the capacitor~\cite{Henriksen}. The relative permittivity of sodium is $8.4$~\cite{Nalwa}. The full 3D capacitor model is illustrated in Fig.~\ref{fig:Capacitor-practical-3D}, left. The domain $F$ is allocated to the domain containing the sodium contamination, surrounded by the pyrex dielectric material. In this experiment, we investigate the simplification error after replacing the sodium domain with pyrex and our domain of interest $S$ is the lower half of the space between the conductor plates. Our QoI represents the capacity in this domain and is interesting for several engineering applications~\cite{Sakalli}. E.g., it is a characteristic that can be used in the design of capacitors or analyze the properties of a molecular system.

\begin{figure}[t]
 \centering
 \includegraphics[width=\textwidth]{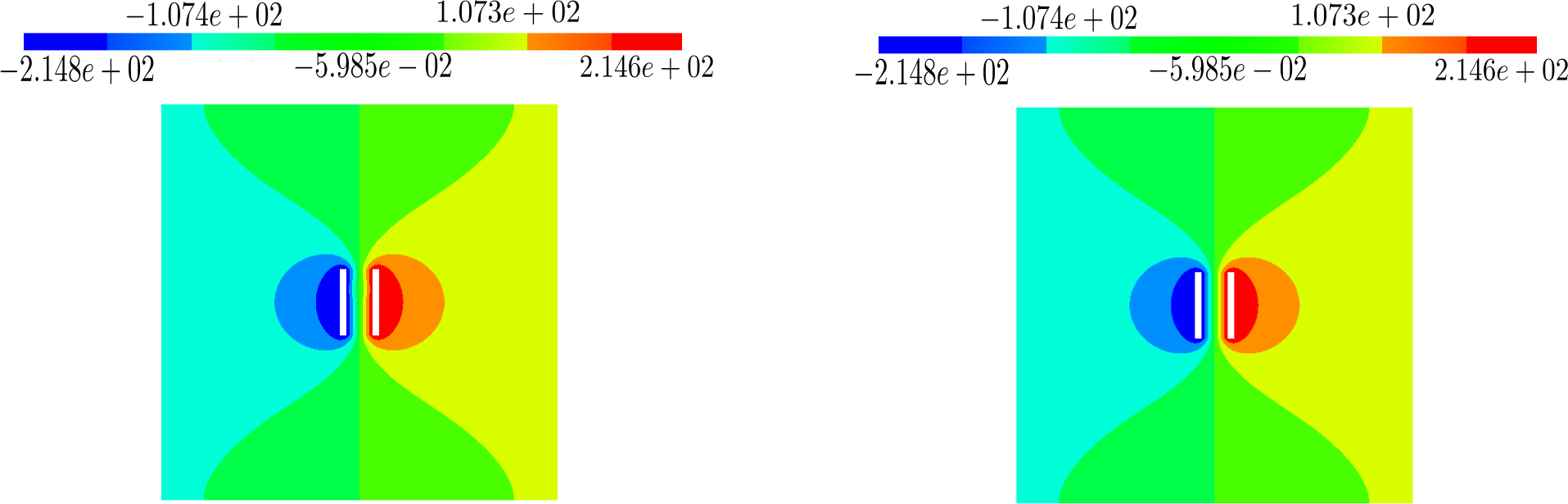}
 \caption{Left: Finite element solution for the original glass capacitor model. The feature affects the distribution of the electrostatic potential nearby. Right: Finite element solution for the capacitor model after defeaturing. The electric field is symmetric between the capacitor plates.}
 \label{fig:Capacitor-Practical-Feature-Defeature}
\end{figure}

The sodium contamination is a cube of length $0.2$cm. The contamination cannot be embedded in or interfaced to the domain of interest $S$ because of the assumptions made for the simplification error estimation. The location of the contamination $F$ in the dielectric material is shown in the Fig.~\ref{fig:Capacitor-practical-3D}(right), showing a cut plane through the middle of the capacitor. It can be seen that the dielectric material is placed between the capacitor plates. The domain of interest $S$ is in the lower half of the  dielectric material where there is no contamination in that region. We run the finite element analysis for both the original and simplified problems. Fig.~\ref{fig:Capacitor-Practical-Feature-Defeature} shows the solution of both. It can be seen that the sodium contamination disturbs the electric field in the pyrex material.

\begin{figure}[t]
 \centering
 \includegraphics[width=\textwidth]{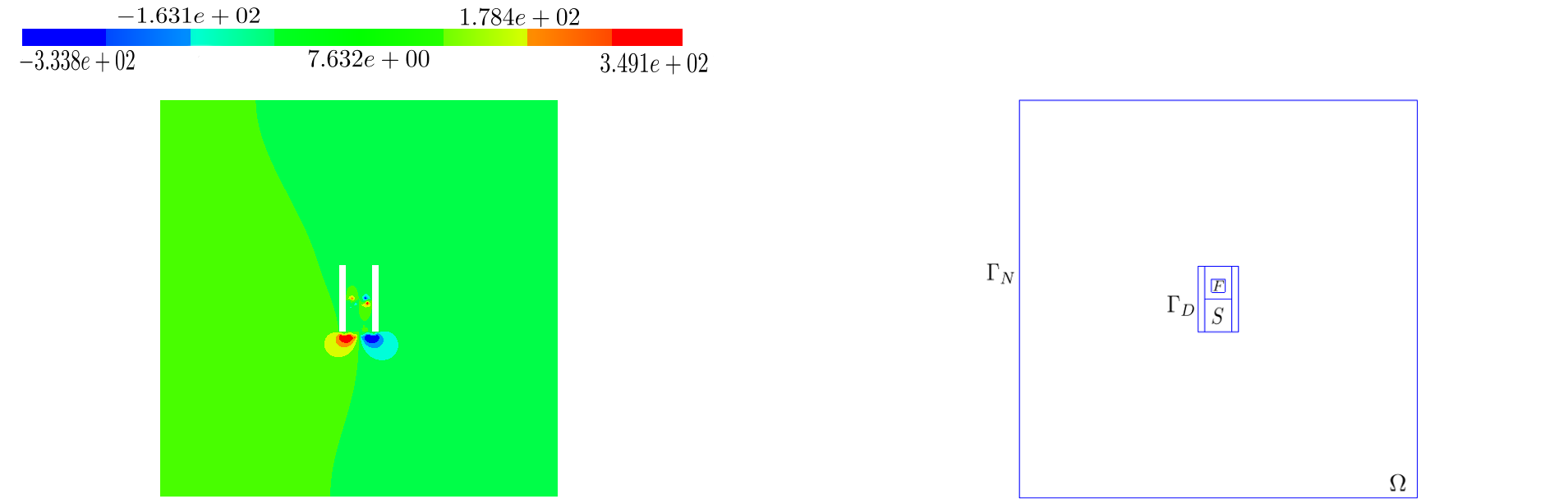}
 \caption{The finite element solution of the dual simplified model for the glass capacitor (left). The distribution of the solution is compared to the referential geometry setup on the right.}
 \label{fig:Capacitor-Practical-Dual-Geo}
\end{figure}

The error bound is computed according to Eq.~\eqref{eq:gb}. It requires the finite element solution of primal and dual simplified models, shown in Figs.~\ref{fig:Capacitor-Practical-Feature-Defeature}, \ref{fig:Capacitor-Practical-Dual-Geo}, to calculate $\nu$, $\nu^\ast$, and $R(\tilde{\Phi}^\ast)$. We are able to obtain the exact and simplified QoI from the solutions of the original and simplified models. The results of this gives $L=293032.95 \leq Q(\Phi)=293213 \leq U=293408.05$ resulting in $\mathfrak{I} = 1.0012793$, indicating good quality bounds. As we calculate the QoI to represent electrostatic energy by approximating Eq.~\eqref{eq:q} using relative permittivities on centimeter scales, these and all following QoI values are in $c * \text{Joules}$ with $c = 2 * 10^4 / \epsilon_0 \approx 2.26 * 10^{15}$ and $\epsilon_0$ is the (dimensionless) electric constant (permittivity of vacuum).

\subsubsection{Experiment 1}

This experiment is used to demonstrate the viability of the method, rather than being a practically useful problem. However, the capacitor model and its condition is similar to the previous model for consistency between examples. We run the experiment for several dielectric materials with different relative permittivities. Here we choose the relative permittivities from $\epsilon_F = 3\epsilon_R$ to $13\epsilon_R$, where $\epsilon_R$ is the relative permittivity of air ($1.0005$). The purpose of this experiment is to test our method to bound the simplification error for different permittivity values in the feature domain where the simplification replaces the material of the feature domain with the relative permittivity of the surrounding domain. We demonstrate how the QoI can be estimated by bounding the simplification error.

\begin{figure}[t]
 \centering
 \includegraphics[width=\textwidth]{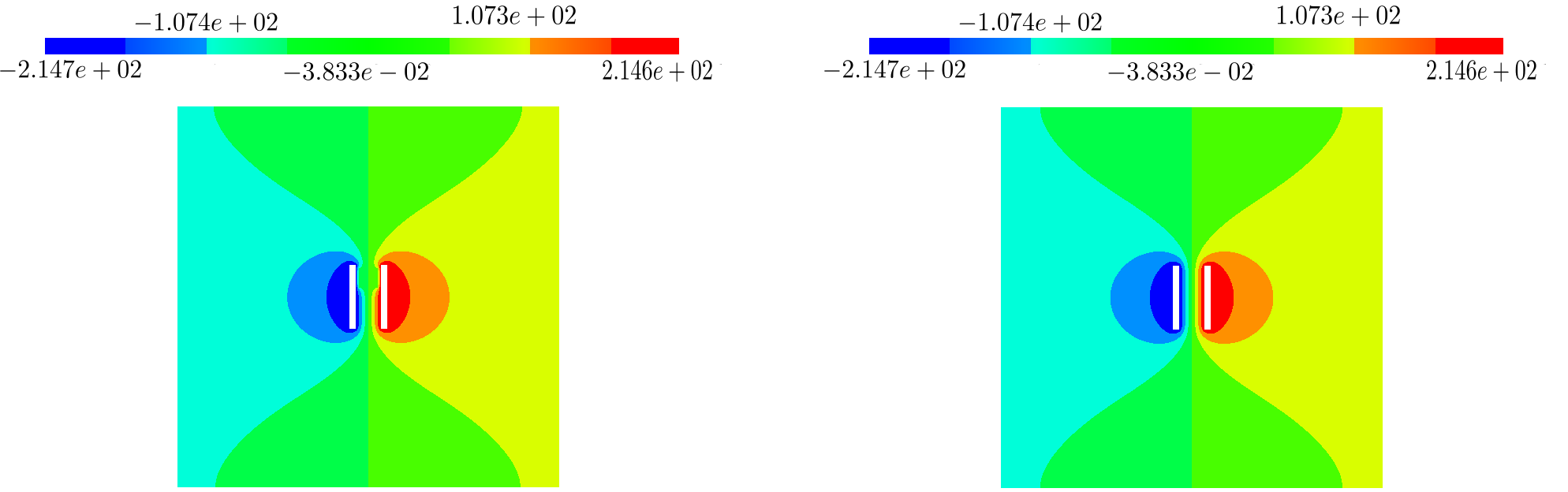}
 \caption{Left: Finite element solution for the original capacitor model for Experiment~1. The feature affects the distribution of the electrostatic potential nearby. Right: Finite element solution for the simplified capacitor model. The field is not distorted in the feature domain and its surrounding domain.}
 \label{fig:Capacitor-defeature-feature-solution}
\end{figure}

\begin{figure}[t]
 \centering
 \includegraphics[width=\textwidth]{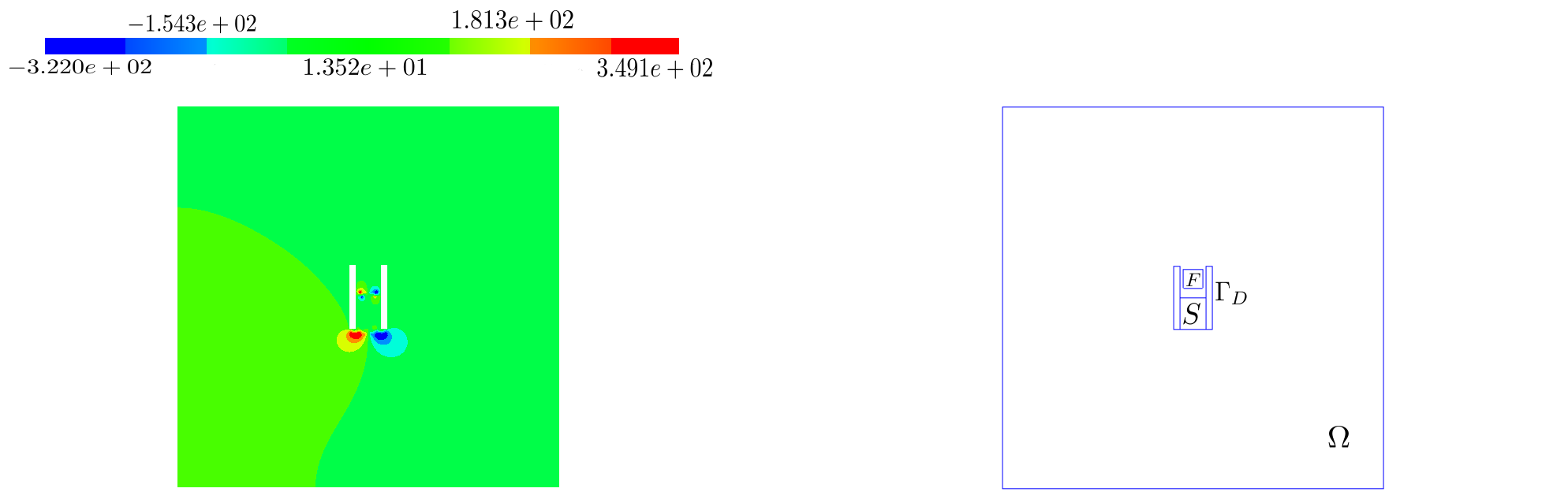}
 \caption{The finite element solution of the adjoint model for Experiment~1 (left) and a representation of the geometry of the model (right).}
 \label{fig:Capacitor-Dual}
\end{figure}

Fig.~\ref{fig:Capacitor-defeature-feature-solution} shows the distribution of the electric potential before and after defeaturing. It can be seen by simple comparison of the two solutions that the dielectric material in the feature domain pushes away the electric field. Removing of dielectric material makes changes in the stored value of electric energy in the domain of interest, $S$. In order to estimate the simplification error we must solve the dual problem represented by Eq.~\eqref{eq:dual}. Its solution is shown in Fig.~\ref{fig:Capacitor-Dual}.

\begin{figure}[t]
  \centering
  \begin{tikzpicture}
  \begin{axis}[width=.45\textwidth,
               height=6cm,
               only marks,
               scaled ticks = true,
               mark size=1pt,
               axis lines = left,
               y axis line style = {blue},
               xlabel = {\small $\epsilon_F$},
               xlabel style = {yshift=5pt},
               xtick distance = 1,
               xticklabel = {{\small$\pgfmathprintnumber{\tick}$}},
               ylabel = {{\small\color{blue}$Q(\Phi)$}},
               ylabel style = {overlay,yshift=-5pt},
               ytick distance = 200000,
               yticklabel = {{\small$\pgfmathprintnumber{\tick}$}},
               every y tick/.style = {blue},
               yticklabel style = {blue},
               yticklabel style= {/pgf/number format/precision=5},
               xmin=2.5,
               xmax=13.5]
  \addplot+[blue,
            error bars/.cd,
            y dir = both, y explicit]
  table[x = e,
        y = q,
        y error plus expr=\thisrow{u}-\thisrow{q},
        y error minus expr=\thisrow{q}-\thisrow{l}] {
    e  q      u       l
    3  243729 322114  204100
    4  243748 372406  173208
    5  243769 424159  140855
    6  243789 476665  107749
    7  243807 529569   74246
    8  243825 582738   40477
    9  243841 613448   29168
    10 243856 689565  -27549
    11 243869 766425 -629861
    12 243882 825171 -723693
    13 243893 884180 -817786
  };
  \end{axis}
  \begin{axis}[red,
               width=.45\textwidth,
               height=6cm,
               only marks,
               scaled ticks = true,
               mark size=1pt,
               axis lines = right,
               axis x line = none,
               xlabel = \empty,
               xtick = \empty,
               ylabel = {{\small$\mathfrak{I}$}},
               ylabel style = {overlay,yshift=5pt},
               yticklabel = {{\small$\pgfmathprintnumber{\tick}$}},
               every y tick/.style = {red},
               ytick distance = 1,
               xmin=2.5,
               xmax=13.5,
               ymin=1,
               ymax=8]
  \addplot+[red, mark=*, mark options = {red}]
  table[x = e,
        y expr = 1 + (\thisrow{u} - \thisrow{l}) / \thisrow{q} ] {
    e  q      u       l
    3  243729 322114  204100
    4  243748 372406  173208
    5  243769 424159  140855
    6  243789 476665  107749
    7  243807 529569   74246
    8  243825 582738   40477
    9  243841 613448   29168
    10 243856 689565  -27549
    11 243869 766425 -629861
    12 243882 825171 -723693
    13 243893 884180 -817786
  };
  \end{axis}
  \end{tikzpicture}
  \hfil\small
  \begin{tabular}[b]{|r|r|r|r|r|}
    \hline
    $\epsilon_F$ & $\mathfrak{I}$ & $Q(\Phi)$ & $U$ & $L$\\ \hline
    \resrow{3}{243729}{322114}{204100}
    \resrow{4}{243748}{372406}{173208}
    \resrow{5}{243769}{424159}{140855}
    \resrow{6}{243789}{476665}{107749}
    \resrow{7}{243807}{529569}{74246}
    \resrow{8}{243825}{582738}{40477}
    \resrow{9}{243841}{613448}{29168}
    \resrow{10}{243856}{689565}{-27549}
    \resrow{11}{243869}{766425}{-629861}
    \resrow{12}{243882}{825171}{-723693}
    \resrow{13}{243893}{884180}{-817786}
    \multicolumn{5}{l}{\vspace*{3ex}}
  \end{tabular}
 \caption{Bounds of the simplification error for the QoI for different dielectric values in an internal feature for Experiment~1. The table lists the relative permittivity $\epsilon_F$ in the feature domain, the effectivity index $\mathfrak{I}$, the exact QoI $Q(\Phi)$, and the upper $U$ and lower $L$ bounds.}
 \label{fig:ErrorEpsilon}
\end{figure}
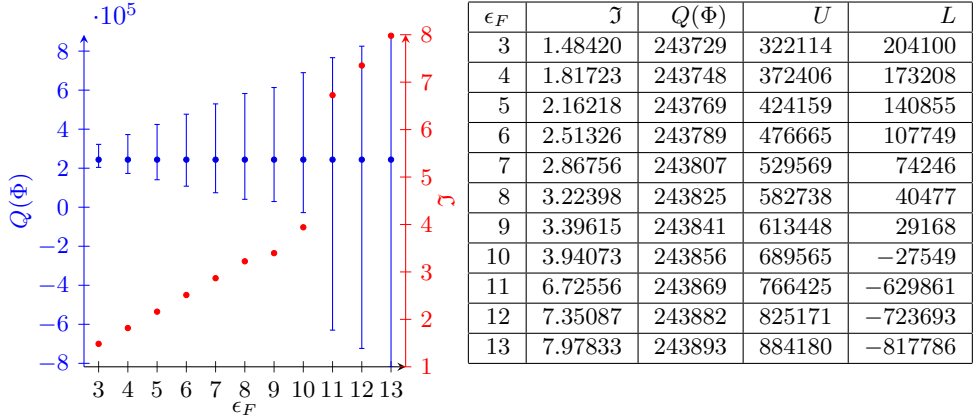

Fig.~\ref{fig:ErrorEpsilon} shows the upper and lower energy norm bounds compared to the exact QoI value. The larger the relative permittivity value, the less tight the estimation of the QoI is. It can even become negative, which is not acceptable. So the simplification error estimation becomes less useful when the influence of the suppressed feature on the QoI becomes stronger. While this is not ideal, the overestimation of the error can still prevent us from making an inappropriate simplification.

\subsubsection{Experiment 2}

\begin{figure}[t]
 \centering
 \includegraphics[width=.35\textwidth]{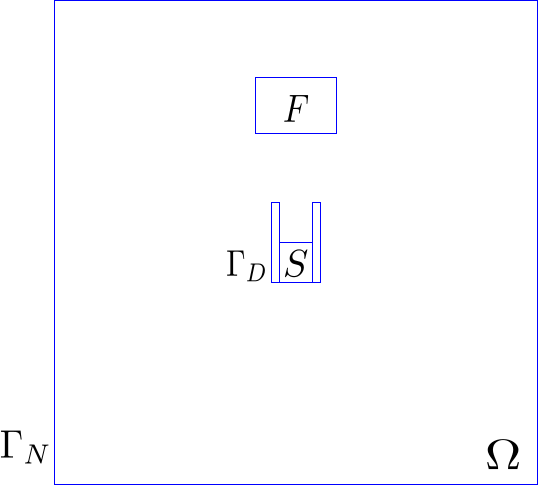}
 \caption{Capacitor with an internal feature of variable size used in Experiment~2.}
 \label{fig:Capacitor-Feature-Size}
\end{figure}

The aim of this experiment is to test the performance of the simplification error for different feature geometries, in particular for a feature increasing in size, under constant relative permittivity in the feature. In this case, the feature $F$ lies outside the conductor plates, while $S$, the capacitor geometry and boundary conditions are the same as in Experiment~1; see Fig.~\ref{fig:Capacitor-Feature-Size}. The relative permittivity of the dielectric material inside $F$ is set to $5\epsilon_R$ for all simulations, where $\epsilon_R$ is the relative permittivity of air ($1.0005$). The capacitor is filled with air. Defeaturing replaces the relative permittivity in the feature domain $F$ with the relative permittivity of the material surrounding it. We increase the size of $F$ from a width and height of $0.5$cm, $0.4$cm to a width and height of $1.6$cm, $1.5$ cm in $0.1$cm steps added to both at the same time. We have also added two larger rectangular feature domains. In each step, the simplification error is bounded as before.

\begin{figure}[t]
 \centering
 \includegraphics[width=\textwidth]{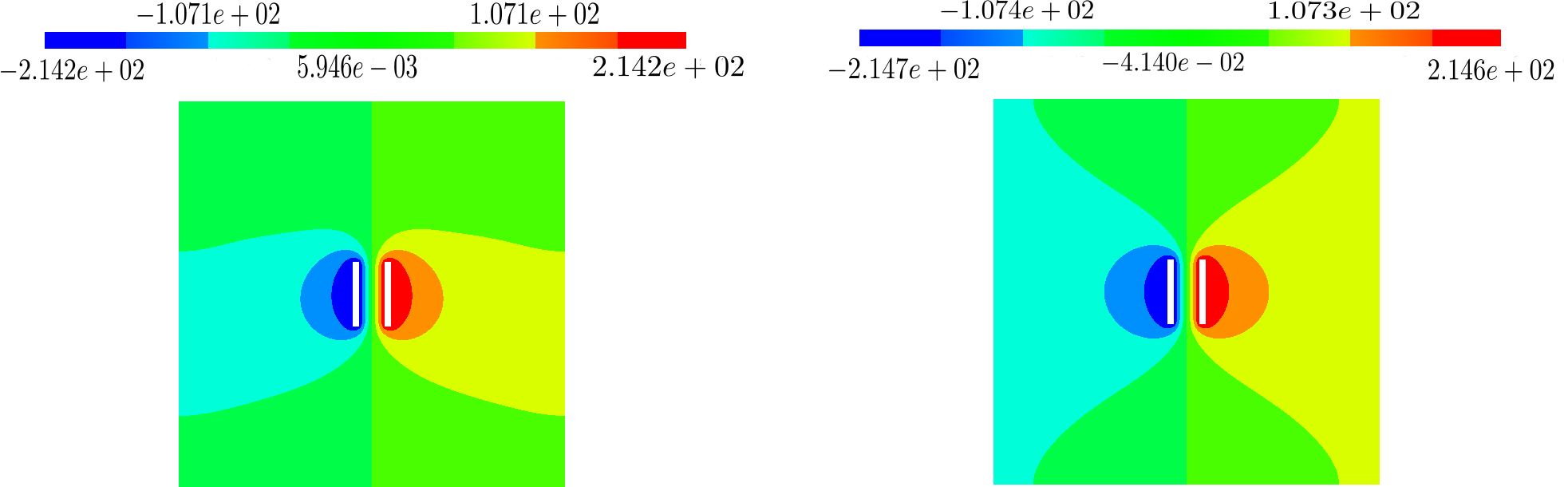}
 \caption{Left: The finite element solution for the field $\Phi$ for the original capacitor model for Experiment~2. Right: The finite element solution of the simplified capacitor model.}
 \label{fig:Capacitor-Feature-Size-Def-Solution}
\end{figure}

As the feature increases in size, it repels the electrostatic potential more, which changes the energy in the domain of interest. The simulation analysis result for the largest feature is shown in Fig.~\ref{fig:Capacitor-Feature-Size-Def-Solution}.

\begin{figure}[t]
  \centering
  \begin{tikzpicture}
  \begin{axis}[width=.33\textwidth,
               height=6cm,
               only marks,
               scaled ticks = true,
               mark size=1pt,
               axis lines = left,
               y axis line style = {blue},
               xlabel = {\small $H * W$},
               xlabel style = {yshift=5pt},
               xtick distance = 1,
               xticklabel = {{\small$\pgfmathprintnumber{\tick}$}},
               xticklabels = {0,0,1,2,4},
               extra x ticks = { 2.5 },
               extra x tick style={grid=none,tick label style={xshift=0cm,yshift=10pt}},
               extra x tick label = {\color{black}{/\!\!/}},
               ylabel = {{\small\color{blue}$Q(\Phi)$}},
               ylabel style = {overlay,yshift=-5pt},
               ytick distance = 40,
               yticklabel = {{\small$\pgfmathprintnumber{\tick}$}},
               every y tick/.style = {blue},
               yticklabel style = {blue},
               yticklabel style= {/pgf/number format/precision=5},
               xmin=0.1,
               xmax=3.3]
  \addplot+[blue,
            error bars/.cd,
            y dir = both, y explicit]
  table[x expr = \thisrow{w} * \thisrow{h},
        y = q,
        y error plus expr=\thisrow{u}-\thisrow{q},
        y error minus expr=\thisrow{q}-\thisrow{l}] {
    w   h    q      u      l
    0.5 0.5  243706 243722 243704
    0.6 0.6  243714 243740 243713
    0.7 0.7  243715 243748 243713
    0.8 0.8  243725 243768 243722
    0.9 0.9  243721 243777 243718
    1.0 1.0  243721 243789 243718
    1.1 1.1  243713 243794 243709
    1.2 1.2  243719 243812 243715
    1.3 1.3  243724 243830 243719
    1.4 1.4  243730 243848 243725
    1.5 1.5  243720 243849 243714
  };
  \addplot+[blue, mark = *, mark options = {blue},
            error bars/.cd,
            y dir = both, y explicit]
  table[x expr = \thisrow{w} * \thisrow{h} - 1,
        y = q,
        y error plus expr=\thisrow{u}-\thisrow{q},
        y error minus expr=\thisrow{q}-\thisrow{l}] {
    w   h    q      u      l
    1.3 2.85 243729 243870 243722
    1.4 2.95 243732 243886 243726
  };
  \end{axis}
  \begin{axis}[red,
               width=.33\textwidth,
               height=6cm,
               only marks,
               scaled ticks = true,
               mark size=1pt,
               axis lines = right,
               axis x line = none,
               xlabel = \empty,
               xtick = \empty,
               ylabel = {{\small$(\mathfrak{I} - 1) * 10^4$}},
               ylabel style = {overlay,yshift=5pt},
               yticklabel = {{\small$\pgfmathprintnumber{\tick}$}},
               yticklabel style= {/pgf/number format/precision=6},
               every y tick/.style = {red},
               ytick distance = 1,
               xmin=0.1,
               xmax=3.3,
               ymin=0,
               ymax=7]
  \addplot+[red, mark=*, mark options = {red}]
  table[x expr = \thisrow{w} * \thisrow{h},
        y expr = 10000 * (\thisrow{u} - \thisrow{l}) / \thisrow{q} ] {
    w   h    q      u      l
    0.5 0.5  243706 243722 243704
    0.6 0.6  243714 243740 243713
    0.7 0.7  243715 243748 243713
    0.8 0.8  243725 243768 243722
    0.9 0.9  243721 243777 243718
    1.0 1.0  243721 243789 243718
    1.1 1.1  243713 243794 243709
    1.2 1.2  243719 243812 243715
    1.3 1.3  243724 243830 243719
    1.4 1.4  243730 243848 243725
    1.5 1.5  243720 243849 243714
  };
  \addplot+[red, mark=*, mark options = {red}]
  table[x expr = \thisrow{w} * \thisrow{h} - 1,
        y expr = 10000 * (\thisrow{u} - \thisrow{l}) / \thisrow{q} ] {
    w   h    q      u      l
    1.3 2.85 243729 243870 243722
    1.4 2.95 243732 243886 243726
  };
  \draw[black] decorate [decoration={zigzag}] {(axis cs:2.5,0) -- (axis cs:2.5,7)};
  \end{axis}
  \end{tikzpicture}
  \hfil\small
  \begin{tabular}[b]{|r|r|r|r|r|}
    \hline
    $(H,W)$ & $\mathfrak{I}$ & $Q(\Phi)$ & $U$ & $L$ \\ \hline
    \resrow{(0.5,0.4)}{243706}{243722}{243704}
    \resrow{(0.6,0.5)}{243714}{243740}{243713}
    \resrow{(0.7,0.6)}{243715}{243748}{243713}
    \resrow{(0.8,0.7)}{243725}{243768}{243722}
    \resrow{(0.9,0.8)}{243721}{243777}{243718}
    \resrow{(1.0,0.9)}{243721}{243789}{243718}
    \resrow{(1.1,1.0)}{243713}{243794}{243709}
    \resrow{(1.2,1.1)}{243719}{243812}{243715}
    \resrow{(1.3,1.2)}{243724}{243830}{243719}
    \resrow{(1.4,1.3)}{243730}{243848}{243725}
    \resrow{(1.5,1.4)}{243720}{243849}{243714}
    \resrow{(1.3,2.85)}{243729}{243870}{243722}
    \resrow{(1.4,2.95)}{243732}{243886}{243726}
  \end{tabular}
  \caption{Bounds on the simplification error in Experiment~2 for varying feature dimensions. The table lists the height $H$ and width $W$ of the feature domain, the effectivity index $\mathfrak{I}$, the exact QoI $Q(\Phi)$, and the upper $U$ and lower $L$ bounds.}
  \label{fig:ErrorExpansion}
\end{figure}
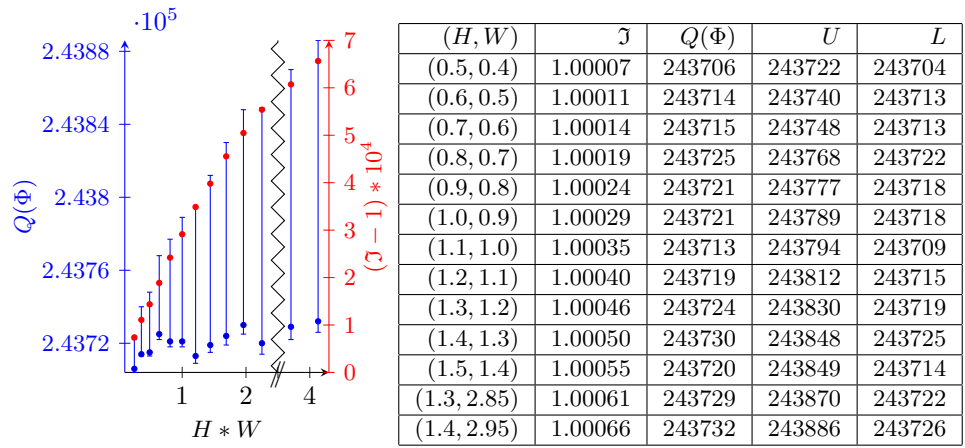

The upper and lower bounds for the QoI with their effectivity indices are shown in the Fig.~\ref{fig:ErrorExpansion}. The effectivity indices in this case are close to one, indicating that the error can be well characterised. The error bounds grow bigger as the size of the feature expands: as the feature $F$ moves further away from the domain of interest $S$, the influence it has on the QoI is reduced. We can conclude that the bounds are tight enough for all the defeaturing experiments and estimate the QoIs well.

Comparing between the effectivity indices for Experiment~1 and Experiment~2, it can be seen that the larger the distance of the feature domain $F$ from the domain of interest $S$, the effect of defeaturing on the QoI is reduced. Increasing the size of the feature reduces the effectivity of the bounds. Similarly, increasing the relative permittivity of the feature, the higher the effect on the QoI, the less effective are the bounds. In all cases the bounds are valid. Albeit, their tightness varies according to the feature condition. The general trend is that the estimation of the error is getting less effective as the effect of the feature on the QoI becomes stronger. This means the bounds are practically useful to make defeaturing decisions.

\section{Boundary Features} \label{Boundary}

In this section we consider features on the boundary of $\Omega$. We distinguish positive features, that are subtracted from $\Omega$ by the simplification, and negative features, that are added to $\Omega$ by the simplification. As in all these cases the simplification implies a change of the domain, we must suitably expand or reduce the defeatured domains $\tilde\Omega$ such that they remain compatible over the original domain, i.e. construct $\widehat{\Omega} = \Omega$, which was not necessary for internal features. We further have to distinguish between features located on the Neumann and Dirichlet boundaries. We first consider the case of Neumann boundary conditions for positive and negative features and then also discuss the Dirichlet boundary condition case.

\subsection{Positive Features with Neumann Boundary Conditions}

We first consider the case of a positive feature being simplified on the Neumann boundary. The domain $\Omega$ is partitioned into a feature domain $F$ and the remaining defeatured domain $\tilde\Omega$ such that $\Omega = \tilde{\Omega} \cup F$. Here defeaturing removes the positive feature by setting the electric permittivity in the feature domain, $\varepsilon_r|_F \equiv 0$, which means the simplified domain $\widehat{\Omega}$ remains equal to $\Omega$, but with different relative permittivity.

\begin{figure}[t]
 \centering
 \includegraphics[width=.35\textwidth]{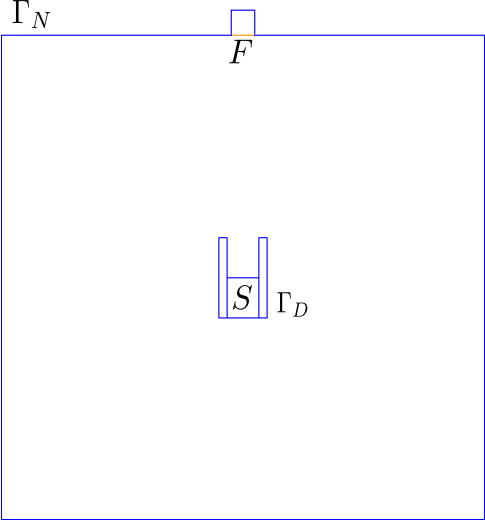}
 \caption{Capacitor model with positive feature $F$ on the Neumann boundary $\Gamma_N$ and domain $\Omega = \tilde{\Omega} \cup F$ for Experiment~3.}
 \label{Capacitor-Pos}
\end{figure}

Note, in order to calculate the error bound, it is required to setup a new boundary value problem on the feature domain and to calculate $\nu$, $\nu^\ast$ and $R(\widehat{\Phi}^\ast)$ from the finite element solution of the model. For this we especially define a Dirichlet boundary condition on the interface between $\tilde{\Omega}$ and $F$, determined by the solution of the defeatured model. So the interface boundary is assigned with the electrostatic potential $\tilde{\Phi}$, obtained from the finite element solution of the defeatured model at the interface between $\tilde{\Omega}$ and $F$, yielding the following problem:
\begin{equation}
\begin{aligned}
 \n \cdot (\epsilon_F \n\Phi_F) &= 0  &&\text{in } F,\\
 \Phi_F &= \tilde{\Phi}  && \text{on } \partial F \cap \partial\tilde{\Omega},\\
 \mathbf{n} \cdot (\epsilon_F \n\Phi_F) &= 0 &&\text{on } \partial F \cap \Gamma_N,
\end{aligned}\label{eq:F}
\end{equation}
where $\Phi_F$ is the solution for this \emph{feature problem}; also see Fig.~\ref{Capacitor-Pos}. At the interface between the feature and defeatured domain, the flux should be equal in length, but point in opposite directions across the interface for the defeatured problem and the feature problem. This means the electrostatic potential remains continuous across the interface for the defeatured problem solution and the feature problem solution which gives an overall simplified solution on $\widehat{\Omega}$. Nevertheless, we assume $D_F = \epsilon_F \n \Phi_F$ is null in $F$, which enables us to calculate the bounds from the simplified model. The rationale for this assumption lies in the construction of the bound based on a posteriori error estimation from the solution of the simplified model and the CRE. To calculate the terms needed for the bound, we require a field over the feature domain and this approach represents effectively the removal of an area from the problem domain. Consequently $\nu$, $\nu^\ast$ and $R(\widehat{\Phi}^\ast)$ are computed over the feature domain where the simplification error bounds are calculated based on the assumption of eliminating the flux. As explained in Section~\ref{General}, the calculation of the error bounds follows Eq.~\eqref{eq:gb}. The calculation of the terms for this feature type is described below.

$\nu$ is related to the CRE Eq.~\eqref{eq:nu} based on the separation of the energy norm computation for $F$ from $\tilde\Omega$. We use the solution $\widehat{\Phi}$ over the simplified domain $\widehat{\Omega}$ which combines the solution of the defeatured problem and the feature problem. Based on Eqs.~\eqref{eq:nu} and~\eqref{eq:nunor}, with the specific relation between $F$ and $\tilde{\Omega}$ in this case, we get
\begin{equation}
\begin{aligned}
\nu^2 =& \| \widehat{D} + \epsilon_R \n \widehat{\Phi} \|^2_{{\epsilon_R}^{-1}} + \| \widehat{D} - \epsilon_F \n \widehat{\Phi} \|^2_{{\epsilon_F}^{-1}} \\
     = &\int_{\tilde{\Omega}} (\widehat{D} - D) {\epsilon_R}^{-1} (\widehat{D} - D) \; \mathrm{d}\Omega\\
        & \qquad {} + \int_F (\widehat{D} - D) {\epsilon_F}^{-1} (\widehat{D} - D) \; \mathrm{d}\Omega + \int_{\tilde{\Omega}} (\n \widehat{\Phi} - \n\Phi) {\epsilon_R}^{-1} (\n \widehat{\Phi} - \n\Phi) \; \mathrm{d}\Omega \\
        & \qquad {} + \int_F (\n \widehat{\Phi} - \n\Phi) {\epsilon_F}^{-1} (\n \widehat{\Phi} - \n\Phi) \; \mathrm{d}\Omega - 2 \underbrace{\int_\Omega (\widehat{D} - D)(\n \widehat{\Phi} - \n\Phi) \;\mathrm{d}\Omega}_{=0}.
\end{aligned} \label{eq:posnu}
\end{equation}
The last term above is zero because of the orthogonality between static and kinematic admissibility conditions. For this to be fulfilled, the following must be true for the flux, though:
\begin{equation}
  \int_\Omega \widehat{D} \cdot \n\Phi \; \mathrm{d}\Omega = 0 \; \forall \Phi \in H^1_0(\Omega). \label{eq:divbo}
\end{equation}
This condition is fulfilled, if $\widehat{D}$ is divergence free and the Neumann boundary condition is fulfilled for the original model. The construction of $\widehat{D}$ means $\widehat{D}|_F \equiv 0$, so it is sufficient that $\widehat{D}|_{\tilde{\Omega}} \equiv \tilde{D}|_{\tilde{\Omega}}$. Hence, we get overall
\begin{equation}
\begin{aligned}
\nu^2 =& \int_{F} (\widehat{D} - \epsilon_F \nabla\widehat{\Phi}) {\epsilon_F}^{-1} (\widehat{D} - \epsilon_F \nabla\widehat{\Phi}) \; \mathrm{d}\Omega + \int_{\widetilde{\Omega}} (\widehat{D} - \epsilon_R \nabla\widehat{\Phi}) {\epsilon_R}^{-1} (\widehat{D} - \epsilon_R \nabla\widehat{\Phi}) \; \mathrm{d}\Omega \\
      =& \int_{F} \epsilon_F \n\widehat{\Phi}\n\widehat{\Phi} \; \mathrm{d}\Omega + \int_{\tilde\Omega} \underbrace{(\widehat{D} - \epsilon_R \nabla\widehat{\Phi}) {\epsilon_R}^{-1} (\widehat{D} - \epsilon_R \n\widehat{\Phi})}_{=0} \; \mathrm{d}\Omega.
\end{aligned}
\end{equation}

$\nu^\ast$ follows the methodology for $\nu$. For the dual problem, $\widehat{\Phi}$ is replaced by $\widehat{\Phi}^\ast$ in Eq.~\eqref{eq:dual}:
\begin{equation}
\begin{aligned}
D^\ast - \tilde{D}_s &= -\varepsilon \nabla\widehat{\Phi}^\ast &&\text{in } S, \\
D^\ast &= -\varepsilon \nabla\widehat{\Phi}^\ast &&\text{in } \widehat{\Omega}\setminus S,\\
\nabla \cdot (- \varepsilon \nabla \widehat{\Phi}^\ast) &= Q(\Psi) &&\text{in } \widehat{\Omega}, \\
\widehat{\Phi}^\ast &= 0 &&\text{on } \Gamma_D,\\
\mathbf{n} \cdot (\varepsilon \nabla\widehat{\Phi}^\ast) &= 0 &&\text{on } \Gamma_N.
\end{aligned}\label{eq:posdual}
\end{equation}
The above dual model has the solution $\widehat{\Phi}^\ast$ from which $\nu^\ast$ is calculated as
\begin{equation}
(\nu^\ast)^2 = \int_{F} \epsilon_F \n\widehat{\Phi}^\ast \n\widehat{\Phi}^\ast \; \mathrm{d}\Omega.
\end{equation}

The residual given by Eq.~\eqref{eq:Res} must be adapted to the positive feature due to the separate feature problem in Eq.~\eqref{eq:F}. It becomes
\begin{equation} \label{eq:boundres}
  R(\widehat{\Phi}^\ast)
   = - \int_{\Omega} \varepsilon_r\n\widehat{\Phi}\n\widehat{\Phi}^\ast \;\mathrm{d}\Omega
   = - \int_{\tilde{\Omega}} \epsilon_R\n\tilde{\Phi}\n\widehat{\Phi}^\ast \;\mathrm{d}\Omega - \int_{F} \epsilon_F\n\widehat{\Phi}\n\widehat{\Phi}^\ast \;\mathrm{d}\Omega
\end{equation}
as $\widehat{\Phi} \equiv \tilde{\Phi}$ in $\tilde{\Omega}$. The bilinear form for the defeatured problem subjected to the prescribed boundary conditions vanishes on $\tilde{\Omega}$, $\int_{\tilde{\Omega}} \epsilon_R\n\tilde{\Phi}\n\widehat{\Phi}^\ast \;\mathrm{d}\Omega = 0$. Hence,
\begin{equation}
  R(\widehat{\Phi}^\ast) = - \int_{F} \epsilon_F\n\widehat{\Phi}\n\widehat{\Phi}^\ast \;\mathrm{d}\Omega.
\end{equation}

The error bounds can be calculated from this according to Eq.~\eqref{eq:gb}. Note that for the calculation of the integrals in this paper we assume the meshes between simplified and original modal remain compatible, which specifically means in this case that they are compatible across the boundary between $F$ and $\tilde\Omega$, to minimise effects from numerical approximations of the integrals.

\subsubsection{Experiment 3}

\begin{figure}[t]
 \centering
 \includegraphics[width=\textwidth]{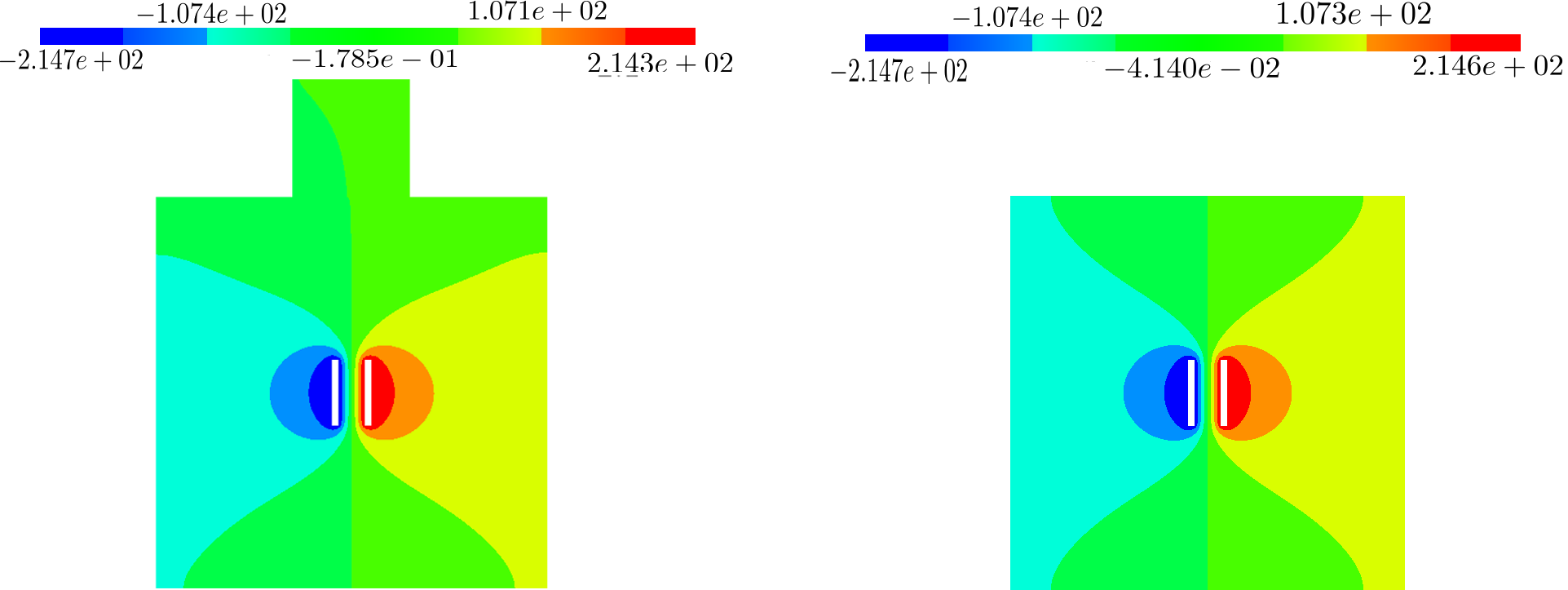}
 \caption{Left: The finite element solution $\Phi$ for the capacitor model with the positive feature on the Neumann boundary for Experiment~3. Right: The finite element solution for the simplified model.}
 \label{fig:Capacitor-Neu-Pos-def-Sol}
\end{figure}

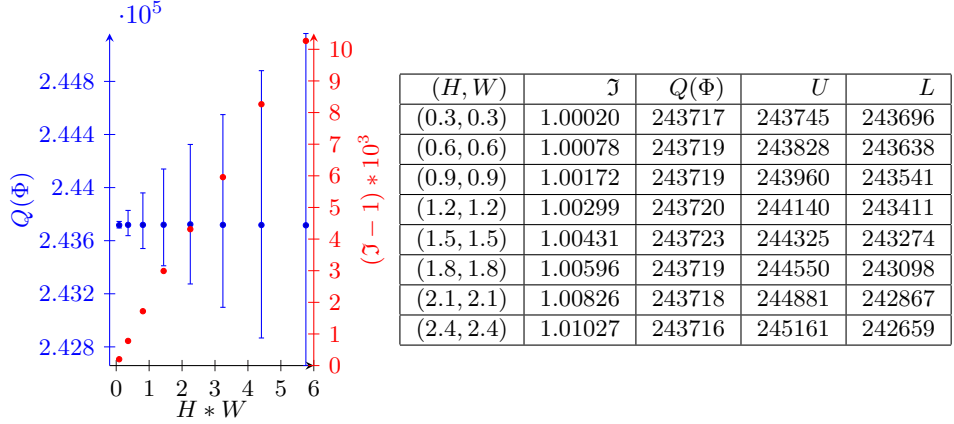
\begin{figure}[t]
  \centering
  \begin{tikzpicture}
  \begin{axis}[width=.33\textwidth,
               height=6cm,
               only marks,
               scaled ticks = true,
               mark size=1pt,
               axis lines = left,
               y axis line style = {blue},
               xlabel = {\small $H * W$},
               xlabel style = {yshift=5pt},
               xtick distance = 1,
               xticklabel = {{\small$\pgfmathprintnumber{\tick}$}},
               ylabel = {{\small\color{blue}$Q(\Phi)$}},
               ylabel style = {overlay,yshift=-5pt},
               ytick distance = 400,
               yticklabel = {{\small$\pgfmathprintnumber{\tick}$}},
               every y tick/.style = {blue},
               yticklabel style = {blue},
               yticklabel style= {/pgf/number format/precision=5},
               xmin=-0.2,
               xmax=6]
  \addplot+[blue,
            error bars/.cd,
            y dir = both, y explicit]
  table[x expr = \thisrow{w} * \thisrow{h},
        y = q,
        y error plus expr=\thisrow{u}-\thisrow{q},
        y error minus expr=\thisrow{q}-\thisrow{l}] {
    w   h   q      u      l
    0.3 0.3 243717 243745 243696
    0.6 0.6 243719 243828 243638
    0.9 0.9 243719 243960 243541
    1.2 1.2 243720 244140 243411
    1.5 1.5 243723 244325 243274
    1.8 1.8 243719 244550 243098
    2.1 2.1 243718 244881 242867
    2.4 2.4 243716 245161 242659
  };
  \end{axis}
  \begin{axis}[red,
               width=.33\textwidth,
               height=6cm,
               only marks,
               scaled ticks = true,
               mark size=1pt,
               axis lines = right,
               axis x line = none,
               xlabel = \empty,
               xtick = \empty,
               ylabel = {{\small$(\mathfrak{I} - 1) * 10^3$}},
               ylabel style = {overlay,yshift=5pt},
               yticklabel = {{\small$\pgfmathprintnumber{\tick}$}},
               yticklabel style= {/pgf/number format/precision=6},
               every y tick/.style = {red},
               ytick distance = 1,
               xmin=-0.2,
               xmax=6,
               ymin=0,
               ymax=10.5]
  \addplot+[red, mark=*, mark options = {red}]
  table[x expr = \thisrow{w} * \thisrow{h},
        y expr = 1000 * (\thisrow{u} - \thisrow{l}) / \thisrow{q} ] {
    w   h   q      u      l
    0.3 0.3 243717 243745 243696
    0.6 0.6 243719 243828 243638
    0.9 0.9 243719 243960 243541
    1.2 1.2 243720 244140 243411
    1.5 1.5 243723 244325 243274
    1.8 1.8 243719 244550 243098
    2.1 2.1 243718 244881 242867
    2.4 2.4 243716 245161 242659
  };
  \end{axis}
  \end{tikzpicture}
  \hfil\small
  \begin{tabular}[b]{|r|r|r|r|r|}
    \hline
    $(H,W)$   & $\mathfrak{I}$ & $Q(\Phi)$ & $U$ & $L$ \\ \hline
    \resrow{(0.3,0.3)}{243717}{243745}{243696}
    \resrow{(0.6,0.6)}{243719}{243828}{243638}
    \resrow{(0.9,0.9)}{243719}{243960}{243541}
    \resrow{(1.2,1.2)}{243720}{244140}{243411}
    \resrow{(1.5,1.5)}{243723}{244325}{243274}
    \resrow{(1.8,1.8)}{243719}{244550}{243098}
    \resrow{(2.1,2.1)}{243718}{244881}{242867}
    \resrow{(2.4,2.4)}{243716}{245161}{242659}
    \multicolumn{5}{l}{\vspace*{5ex}}
   \end{tabular}
  \caption{Bounds on the simplification error in Experiment~3 for varying boundary feature dimensions. The table lists the height $H$ and width $W$ of the feature domain, the effectivity index $\mathfrak{I}$, the exact QoI $Q(\Phi)$, and the upper $U$ and lower $L$ bounds.}\label{ErrPosNeu}
\end{figure}

To test the bound for a positive feature on the Neumann boundary, we again use the QoI from Section~\ref{QoI} for a capacitor model; see Fig~\ref{Capacitor-Pos}. The positive feature is a square of length $0.3$cm initially. The length is increased by $0.3$cm up to $2.4$cm. The solution of electrostatic potential for the original and simplified model is shown in Fig.~\ref{fig:Capacitor-Neu-Pos-def-Sol}. The results for the simplification error bounds are shown in Fig.~\ref{ErrPosNeu}. Overall the effectivity index still indicates good performance and is in line with the conclusions for the internal feature results.

\subsection{Negative Feature on Neumann Boundary} \label{Neumann.Negative}

\begin{figure}[t]
 \centering
  \includegraphics[width=.35\textwidth]{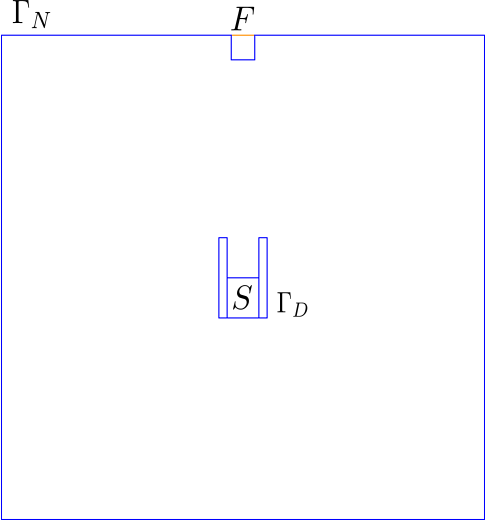}
 \caption{The geometry configuration of the capacitor model with the negative feature $F$ on the Neumann boundary $\Gamma_N$ for Experiment~4.}
 \label{Capacitor-Neg}
\end{figure}

In the case of a negative feature on the Neumann boundary, the simplified domain $\tilde\Omega = \Omega \cup F$ contains $F$, while  $F$ is not a subset of $\Omega$. We must extend $\varepsilon_r$ to $\widehat\Omega = \tilde\Omega = \Omega \cup F$ and set the relative permittivity $\epsilon_F$ of the feature domain close to $0$; see Fig.~\ref{Capacitor-Neg}. With this we make the negative feature part of the domain $\widehat{\Omega}$ and set a permittivity that indicates the original void. This is related to using a penalty factor for the relative permittivity in feature domain,
\begin{equation} \label{eq:ep}
  \widehat{\varepsilon}_r(x) = \alpha \cdot \epsilon_R \text{ for } x \in F,
\end{equation}
where we have chosen $\alpha = 10^{-5}$ for the numerical results and $\epsilon_R$ is the permittivity of the surrounding domain.

We follow the same approach as Section~\ref{General} to bound the simplification error with Eq.~\eqref{eq:gb}. The only difference is that we set a small relative permittivity on $F$. According to Eqs.~\eqref{eq:nu} and~\eqref{eq:nunor} this gives:
\begin{equation}
\begin{aligned}
\nu^2 &= \| \widehat{D} - \epsilon_R \n \widehat{\Phi} \|^2_{{\epsilon_R}^{-1}} + \| \widehat{D} - \alpha\epsilon_R \n \widehat{\Phi} \|^2_{{\epsilon_F}^{-1}} \\
      &= \int_{\tilde{\Omega}} (\widehat{D} - \epsilon_R \n \widehat{\Phi}) {\epsilon_R}^{-1} (\widehat{D} - \epsilon_R \n \widehat{\Phi}) \; \mathrm{d}\Omega\\
      &\qquad {} + \int_F (\widehat{D} - \alpha\epsilon_R \n \widehat{\Phi}) {(\alpha\epsilon_R)}^{-1} (\widehat{D} - \alpha\epsilon_R \n \widehat{\Phi}) \; \mathrm{d}\Omega \\
      &= \int_{F} \alpha \epsilon_R \n\widehat{\Phi}\n\widehat{\Phi} \; \mathrm{d}\Omega,
\end{aligned}
\end{equation}
and similarly, based on Eq.~\eqref{eq:erdua},
\begin{equation}
  (\nu^\ast)^2 = \int_{F} \alpha \epsilon_R \n\widehat{\Phi}^\ast\n\widehat{\Phi}^\ast \; \mathrm{d}\Omega.
\end{equation}
The residual in Eq.~\eqref{eq:residual} is derived similarly to Eq.~\eqref{eq:boundres}, yielding
\begin{equation}\label{eq:resnegbo}
R(\widehat{\Phi}^\ast) = - \int_{F} \alpha \epsilon_R \n\widehat{\Phi}\n\widehat{\Phi}^\ast \;\mathrm{d}\Omega.
\end{equation}

\subsubsection{Experiment 4}

\begin{figure}[t]
 \centering
 \includegraphics[width=\textwidth]{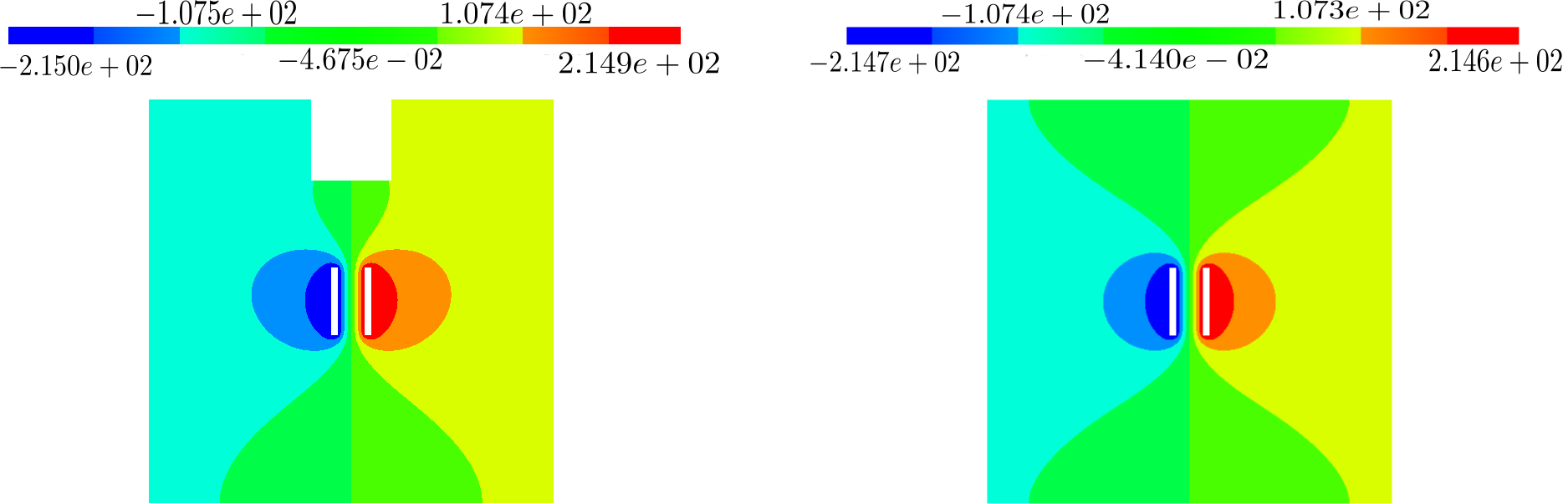}
 \caption{Left: The finite element solution for the capacitor model with a negative feature on the Neumann boundary for Experiment~4. Right: The finite element solution for the simplified model.}
 \label{fig:Capacitor-Neu-Neg-Def-Sol}
\end{figure}

\begin{figure}[t]
  \centering
  \begin{tikzpicture}
  \begin{axis}[width=.33\textwidth,
               height=6cm,
               only marks,
               scaled ticks = true,
               mark size=1pt,
               axis lines = left,
               y axis line style = {blue},
               xlabel = {\small $H * W$},
               xlabel style = {yshift=5pt},
               xtick distance = 1,
               xticklabel = {{\small$\pgfmathprintnumber{\tick}$}},
               ylabel = {{\small\color{blue}$Q(\Phi)$}},
               ylabel style = {overlay,yshift=-5pt},
               ytick distance = 250000,
               yticklabel = {{\small$\pgfmathprintnumber{\tick}$}},
               every y tick/.style = {blue},
               yticklabel style = {blue},
               yticklabel style= {/pgf/number format/precision=5},
               xmin=-0.2,
               xmax=6]
  \addplot+[blue,
            error bars/.cd,
            y dir = both, y explicit]
  table[x expr = \thisrow{w} * \thisrow{h},
        y = q,
        y error plus expr=\thisrow{u}-\thisrow{q},
        y error minus expr=\thisrow{q}-\thisrow{l}] {
    w   h   q      u       l
    0.3 0.3 243719 247400  240034
    0.6 0.6 243715 259020  228414
    0.9 0.9 243713 280545  206883
    1.2 1.2 243709 316112  171281
    1.5 1.5 243702 376369  111031
    1.8 1.8 243693 485712  1686
    2.1 2.1 243681 713065  -225697
    2.4 2.4 243680 1402538 -915189
  };
  \end{axis}
  \begin{axis}[red,
               width=.33\textwidth,
               height=6cm,
               only marks,
               scaled ticks = true,
               mark size=1pt,
               axis lines = right,
               axis x line = none,
               xlabel = \empty,
               xtick = \empty,
               ylabel = {{\small$\mathfrak{I}$}},
               ylabel style = {overlay,yshift=5pt},
               yticklabel = {{\small$\pgfmathprintnumber{\tick}$}},
               yticklabel style= {/pgf/number format/precision=6},
               every y tick/.style = {red},
               ytick distance = 1,
               xmin=-0.2,
               xmax=6,
               ymin=0,
               ymax=11]
  \addplot+[red, mark=*, mark options = {red}]
  table[x expr = \thisrow{w} * \thisrow{h},
        y expr = 1 + (\thisrow{u} - \thisrow{l}) / \thisrow{q} ] {
    w   h   q      u       l
    0.3 0.3 243719 247400  240034
    0.6 0.6 243715 259020  228414
    0.9 0.9 243713 280545  206883
    1.2 1.2 243709 316112  171281
    1.5 1.5 243702 376369  111031
    1.8 1.8 243693 485712  1686
    2.1 2.1 243681 713065  -225697
    2.4 2.4 243680 1402538 -915189
  };
  \end{axis}
  \end{tikzpicture}
  \hfil\small
  \begin{tabular}[b]{|r|r|r|r|r|}
    \hline
    ($W$,$H$)   & $\mathfrak{I}$ & $Q(\Phi)$ & $U$ & $L$ \\ \hline
    \resrow{(0.3,0.3)}{243719}{247400}{240034}
    \resrow{(0.6,0.6)}{243715}{259020}{228414}
    \resrow{(0.9,0.9)}{243713}{280545}{206883}
    \resrow{(1.2,1.2)}{243709}{316112}{171281}
    \resrow{(1.5,1.5)}{243702}{376369}{111031}
    \resrow{(1.8,1.8)}{243693}{485712}{1686}
    \resrow{(2.1,2.1)}{243681}{713065}{-225697}
    \resrow{(2.4,2.4)}{243680}{1402538}{-915189}
    \multicolumn{5}{l}{\vspace*{5ex}}
   \end{tabular}
   \caption{Bounds on the simplification error in Experiment~4 for varying boundary feature dimensions. The table lists the height $H$ and width $W$ of the feature domain, the effectivity index $\mathfrak{I}$, the exact QoI $Q(\Phi)$, and the upper $U$ and lower $L$ bounds.}\label{fig:ErrNegNeu}
\end{figure}
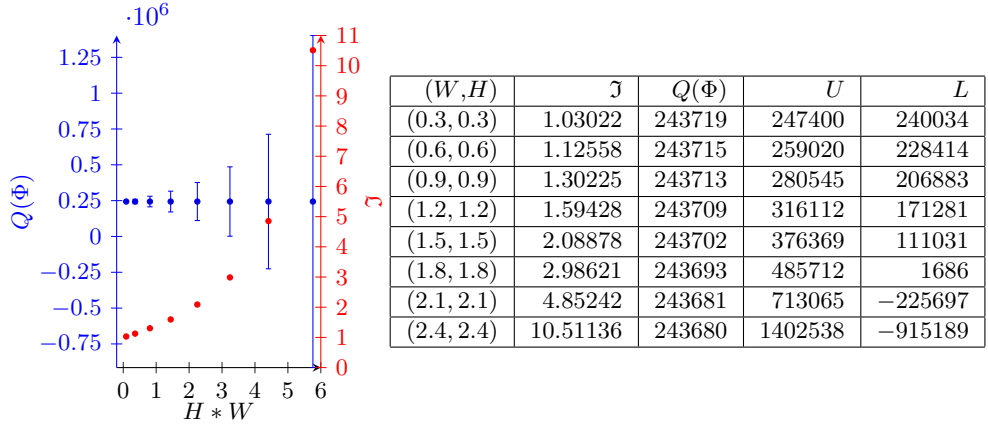

We now test our error bounds for a negative feature on the Neumann boundary with increasing size, with the same QoI and a similar capacitor model as before; see Fig.~\ref{Capacitor-Neg}. The length of the initial square feature is $0.3$cm, increased in steps of $0.3$cm to $2.4$cm. Fig.~\ref{fig:Capacitor-Neu-Neg-Def-Sol} shows the finite element solutions for the simplified and original model. Fig.~\ref{fig:ErrNegNeu} shows the results. We see that by making the slot size bigger, the value of the QoI of the simplified model deviates from the exact one. The values of the effectivity indices also clearly indicate that increasing the size of the negative feature drastically reduces the tightness of the bound. The effectivity index for the largest two features indicate that our simplification error estimation is not suitable for such large features. However, the effectivity indices for smaller sizes feature are in the acceptable range to bound the error. It appears that if the size of the negative feature exceeds a specific threshold, the bounds very quickly become less tight.

\subsection{Negative Feature on Dirichlet Boundary}

\begin{figure}[t]
 \centering
 \includegraphics[width=.16\textwidth]{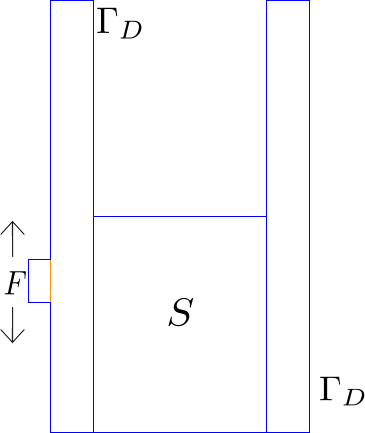}\hfil\includegraphics[width=.2\textwidth]{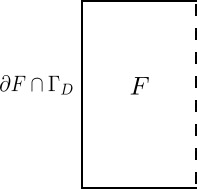}
 \caption{Left: a negative feature on a capacitor plate (Dirichlet Boundary) and its direction of motion for Experiment~5. Right: the standalone negative feature and its associated boundary conditions.}
 \label{fig:Capacitor-Neg-Dir}
\end{figure}

In the capacitor model, the conductors are described by Dirichlet boundary conditions. So for simplifying features on the capacitor plates, we have to consider these cases. Again, the feature can be positive or negative, but as the boundary condition type is different, the general equations in Section~\ref{General} must be adapted differently from the Neumann boundary cases. Here we deal with the simplification error by removing the negative feature on Dirichlet boundary condition only; the positive feature case is similar. In either case, the following main conditions must be fulfilled:
\begin{enumerate}
\item The electrostatic potential (field) must be continuous across the interface between the feature and the rest of the domain.
\item The electrostatic displacement (flux) must be divergence free.
\end{enumerate}
For an example of a negative feature on a capacitor plate see Fig.~\ref{fig:Capacitor-Neg-Dir} (left). The position of the feature is the dominant factor determining the simplification error. As for the negative feature on the Neumann boundary, we must add the negative feature to the domain and with this actually mesh the feature domain. Again, we introduce a penalty factor as in Eq.~\ref{eq:ep} to extend the relative permittivity to the feature domain. The assumption for calculating the energy norms on the feature domain is that all vertices in feature domain are assigned with the Dirichlet electrostatic potential values $\Phi_D$, i.e. the whole feature domain contains a uniform electrostatic potential ($\Phi_F = \Phi_D$). The simplified domain is $\widehat{\Omega} = \tilde{\Omega} \setminus F = \Omega$. The breakup of the simplified domain ensures all integrals are well defined.

The electric displacement (flux) should be divergence free in the feature domain and equal to the flux in the capacitor domain at the interface, though in opposite directions. So overall we have
\begin{align}
  \widehat{\Omega} &= \widetilde{\Omega} \setminus F,\\
  \widehat{\Phi} &= \Phi_D &&\text{in } F, \\
  \widehat{D} &= \tilde{D} &&\text{in } \partial F \cap \Gamma_D.
\end{align}
The derivation of the error bounds is otherwise equivalent to the negative boundary feature case in Section~\ref{Neumann.Negative}, following Eq.~\eqref{eq:gb}. The new assumption and conditions change the results of the integrations, still restricted to the feature domain $F$. The relative permittivity $\epsilon_F$ of the feature domain is adjusted to simulate the closest conditions with respect to the original model. Fig.~\ref{fig:Capacitor-Neg-Dir}(right) depicts the boundary of the feature and its intersection with the Dirichlet boundary condition.

Hence, from Eqs.~\eqref{eq:nu} and~\eqref{eq:nunor} and we get
\begin{equation}
  \nu^2 = \int_{F} \alpha \epsilon_R \n\widehat{\Phi}\n\widehat{\Phi} \; \mathrm{d}\Omega
\end{equation}
under consideration of the field being equal to Dirichlet value in the feature domain. The calculation of $\nu^\ast$ requires the solution $\widehat{\Phi}^\ast$ of the simplified dual model~\eqref{eq:posdual},
\begin{equation}
  (\nu^\ast)^2 = \int_{F} \alpha \epsilon_R \n\widehat{\Phi}^\ast \n\widehat{\Phi}^\ast \; \mathrm{d}\Omega.
\end{equation}
The residual is given by Eq.~\eqref{eq:Res} adapted to the negative feature on the Dirichlet boundary and its new conditions. It gives the same result as Eq.~\eqref{eq:resnegbo}, except that the field in the feature domain is replaced with the Dirichlet boundary value, $\Phi_F = \Phi_D$:
\begin{equation}
  R(\widehat{\Phi}^\ast) = - \int_{F} \alpha \epsilon_R \n\widehat{\Phi}\n\widehat{\Phi}^\ast \;\mathrm{d}\Omega.
\end{equation}

\subsubsection{Experiment 5}

\begin{figure}[t]
 \centering
 \includegraphics[width=\linewidth]{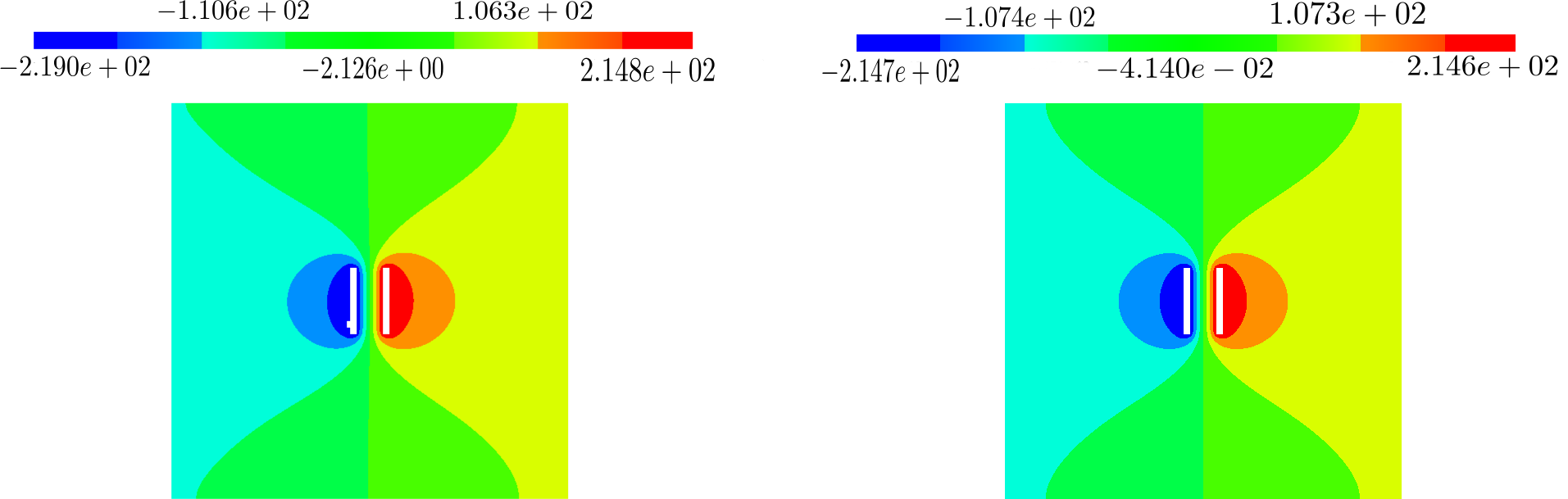}
 \caption{Left: Finite element solution for the capacitor model with a negative feature on the conductor plate (Dirichlet Boundary) for Experiment~5. Right: Finite element solution of the simplified model.}
 \label{fig:Capacitor-Neg-Dir-Def-Sol}
\end{figure}

As shown in Fig.~\ref{fig:Capacitor-Neg-Dir}, the feature moves along the outer side of the left conductor blade. It takes the shape of a negative feature on the conductor which is described by a Dirichlet boundary condition. With the same choice of QoI as before, we are interested in the effect of removing the negative feature, e.g. representing a manufacturing deficiency at different locations, and specifically the performance of our error bounds. The feature is moved along the conductor plate by a step size of $0.1$cm. Fig.~\ref{fig:Capacitor-Neg-Dir-Def-Sol} shows the solution for the simplified and original problem, squeezed in due to the feature.

\begin{figure}[t]
  \centering
  \begin{tikzpicture}
  \begin{axis}[width=.31\textwidth,
               height=6cm,
               only marks,
               scaled ticks = true,
               mark size=1pt,
               axis lines = left,
               y axis line style = {blue},
               xlabel = {\small $Y$},
               xlabel style = {yshift=5pt},
               xtick distance = 0.2,
               xticklabel = {{\small$\pgfmathprintnumber{\tick}$}},
               ylabel = {{\small\color{blue}$Q(\Phi)$}},
               ylabel style = {overlay,yshift=-5pt},
               ytick distance = 25000,
               yticklabel = {{\small$\pgfmathprintnumber{\tick}$}},
               every y tick/.style = {blue},
               yticklabel style = {blue},
               yticklabel style= {/pgf/number format/precision=5},
               xmin=-0.4,
               xmax=0.4]
  \addplot+[blue,
            error bars/.cd,
            y dir = both, y explicit]
  table[x = Y,
        y = q,
        y error plus expr=\thisrow{u}-\thisrow{q},
        y error minus expr=\thisrow{q}-\thisrow{l}] {
    Y     q      u       l
    -0.35 243723 349083 138385
    -0.25 243713 290430 197012
    -0.15 243716 271806 215635
    -0.05 243711 249765 237671
    0.05  243708 254421 233004
    0.15  243715 246556 240884
    0.25  243716 245746 241688
    0.35  243718 245085 242349
  };
  \end{axis}
  \begin{axis}[red,
               width=.31\textwidth,
               height=6cm,
               only marks,
               scaled ticks = true,
               mark size=1pt,
               axis lines = right,
               axis x line = none,
               xlabel = \empty,
               xtick = \empty,
               ylabel = {{\small$(\mathfrak{I}-1)*10$}},
               ylabel style = {overlay,yshift=5pt},
               yticklabel = {{\small$\pgfmathprintnumber{\tick}$}},
               yticklabel style= {/pgf/number format/precision=6},
               every y tick/.style = {red},
               ytick distance = 1,
               xmin=-0.4,
               xmax=0.4,
               ymin=0,
               ymax=9]
  \addplot+[red, mark=*, mark options = {red}]
  table[x = Y,
        y expr = 10 * (\thisrow{u} - \thisrow{l}) / \thisrow{q} ] {
    Y     q      u       l
    -0.35 243723 349083 138385
    -0.25 243713 290430 197012
    -0.15 243716 271806 215635
    -0.05 243711 249765 237671
    0.05  243708 254421 233004
    0.15  243716 246556 240884
    0.25  243716 245746 241688
    0.35  243718 245085 242349
  };
  \end{axis}
  \end{tikzpicture}
  \hfil\small
  \begin{tabular}[b]{|r|r|r|r|r|} \hline
    $(X,Y)$  & $\mathfrak{I}$ & $Q(\Phi)$ & $U$ & $L$ \\ \hline
    \resrow{(-0.325,-0.35)}{243723}{349083}{138385}
    \resrow{(-0.325,-0.25)}{243713}{290430}{197012}
    \resrow{(-0.325,-0.15)}{243716}{271806}{215635}
    \resrow{(-0.325,-0.05)}{243711}{249765}{237671}
    \resrow{(-0.325,0.05)}{243708}{254421}{33004}
    \resrow{(-0.325,0.15)}{243716}{246556}{240884}
    \resrow{(-0.325,0.25)}{243716}{245746}{241688}
    \resrow{(-0.325,0.35)}{243718}{245085}{242349}
    \multicolumn{5}{l}{\vspace*{5ex}}
  \end{tabular}
  \caption{Bounds on the simplification error in Experiment~5 for varying boundary feature location. The table lists the feature location $(X,Y)$, the effectivity index $\mathfrak{I}$, the exact QoI $Q(\Phi)$, and the upper $U$ and lower $L$ bounds.}\label{fig:ErrNegDir}
\end{figure}
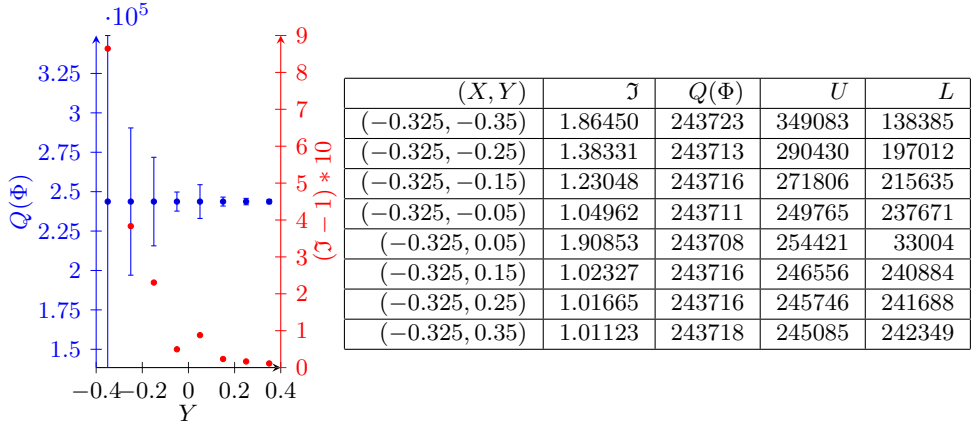

Results for our error bounds, shown in Fig.~\ref{fig:ErrNegDir}, emphasize that the position of the feature plays a significant role in how the simplification affects the QoI. As the feature moves upwards from the bottom of the conductor plate to middle of the plate, where it is next to the domain of interest $S$, the simplification error spikes. After it passed that location, the simplification error and the bound widths are significantly reduced. The effectivity indices indicate that when the feature is next to the domain of interest $S$, the estimation of the simplification error is worse. The error bound effectivity becomes better when the feature moves further upwards from $S$. It also shows that when the feature is at the bottom of the conductor plate, next to $S$, the bound is less tight. This is still in line with our general estimation results earlier that small effects of the feature are quite well estimated, but larger effects are overestimated. Hence, the bounds are a useful practical indicator of when a feature has a small effect on the QoI, but not vice versa.

\section{Conclusion and Future Work} \label{conclusion}

We have devised a novel, conceptually explicit strategy to build lower and upper bounds for the simplification error for divergence free second order Laplacian PDEs. The bounds are based upon the constitutive relation equation which implies the wide applicability of the bounds. We have implemented the bounds to show their effectivity numerically for electrostatics problems, in particular capacitor models. The results are promising: the performance of the bounds for several defeaturing types is good to determine whether the geometry simplification of the model has got only a small effect on a quantity of interest, defined by the user. For larger effects of the geometry simplification on the quantity of interest, the effectivity of the bounds is considerably less tight. However, even if the effect is overestimated in these cases, the bounds are still practically useful to decide whether the effect is small and hence the simplification can be applied.
The requirement that the QoI is linear in the solution of the PDE limits the analysis cases for which the approach can be applied. If the actual QoI can be approximated with a linear QoI such that the approximation error is negligible compared to the simplification error the approach works well, as demonstrated in our numerical experiments. Otherwise, it has to be adapted to non-linear QoIs based on, e.g., non-linear goal-oriented error estimation approaches~\cite{Bryant}. The general strategy presented can be applied to a much wider range of electrostatics and similar problems where the QoIs can be expressed in terms of the energy norm. In future work the effectiveness of our approach for divergence free Laplacian PDEs arising in other electromagnetic or elasticity problems can be explored. QoIs that cannot be expressed in terms of the energy norm, such as eigenmodes, require a different approach, even if the energy norm can still serve as an indication of the error in the field caused by the simplification. Furthermore, the usefulness of the proposed strategy to study the effect of interacting features and their removal is highly relevant.

\end{document}